 \documentclass[final,5p,times,twocolumn]{elsarticle}

\usepackage{amssymb}
\usepackage{cite}
\usepackage{amsmath,amssymb,amsfonts}
\usepackage{algorithmic}
\usepackage{graphicx}
\usepackage{textcomp}
\usepackage{xcolor}
\usepackage{enumitem}
\usepackage[numbers]{natbib}
\usepackage{booktabs}
\usepackage{float}
\usepackage{hyperref}
\usepackage{natbib}
\usepackage{arydshln}
\usepackage{caption}
\usepackage{subcaption}
\usepackage{etoolbox}
\usepackage{comment}
\usepackage{arydshln}

\newcommand{\Real}{\mathbb{R}}
\newcommand{\dt}{\partial t}
\newcommand{\Ho}{\mathcal{H}}
\newcommand{\Bo}{\mathcal{B}}
\newcommand{\Co}{\mathcal{C}}

\newtheorem{theo}{Theorem}
\newtheorem{proposition}{Proposition}
\newtheorem{definition}{Definition}
\newtheorem{assumption}{Assumption}
\newtheorem{remark}{Remark}
\newtheorem{lemma}[theo]{Lemma}
\newenvironment{Proof}[1][Proof]{\textbf{#1.} }{\ \rule{0.5em}{0.5em}}

\newtheorem{example}{Example}
\AtBeginEnvironment{example}{\small}

\definecolor{BlueMatlab}{rgb}{0, 0.4470, 0.7410}
\definecolor{OrangeMatlab}{rgb}{0.8500, 0.3250, 0.0980}
\definecolor{VioletMatlab}{rgb}{0.4940, 0.1840, 0.5560}
\definecolor{GreenMatlab}{rgb}{0.4660, 0.6740, 0.1880}
\definecolor{CyanMatlab}{rgb}{0.3010, 0.7450, 0.9330}

\newcommand{\blue}[1]{\textcolor{BlueMatlab}{#1}}

\newcommand{\violet}[1]{\textcolor{VioletMatlab}{#1}}

\definecolor{DM}{rgb}{0, 0, 0}
\newcommand{\DM}[1]{\textcolor{DM}{#1}}

\journal{Systems $\&$ Control Letters}

\begin{document}

\begin{frontmatter}



\title{Structure-Preserving Discretization and Model Order Reduction of Boundary-Controlled 1D Port-Hamiltonian Systems}

 \author[label1]{Jesus-Pablo Toledo-Zucco\corref{cor1}}
 \cortext[cor1]{Corresponding author} 
\ead[cor1]{jtoledoz@onera.fr}
\ead[url]{https://sites.google.com/view/jesustoledo/main}

 \author[label2]{Denis Matignon}
  \author[label1]{Charles Poussot-Vassal}
 \author[label3]{Yann Le Gorrec}

 \affiliation[label1]{organization={ONERA-DTIS and Fédération ENAC ISAE-SUPAERO ONERA, Université de Toulouse},
             city={Toulouse},
             country={France}}
             
 \affiliation[label2]{organization={Fédération ENAC ISAE-SUPAERO ONERA, Université de Toulouse},
             city={Toulouse},
             country={France}}
             
  \affiliation[label3]{organization={SUPMICROTECH, CNRS, FEMTO-ST},
             city={Besançon},
             postcode={25000},
             country={France}}

%\affiliation{organization={},%Department and Organization

\begin{abstract}
This paper presents a systematic methodology for the discretization and reduction of a class of one-dimensional Partial Differential Equations (PDEs) with inputs and outputs \DM{collocated} at the spatial boundaries. The class of system that we consider is known as Boundary-Controlled Port-Hamiltonian Systems (BC-PHSs) and covers a wide class of Hyperbolic PDEs with a large type of boundary inputs and outputs. This is, for instance, the case of waves and beams with Neumann, Dirichlet, or mixed boundary conditions. Based on a Partitioned Finite Element Method (PFEM), we \DM{develop} a numerical scheme for the structure-preserving spatial discretization for the class of one-dimensional BC-PHSs. We show that if the initial PDE is {\it passive} (or {\it impedance energy preserving}), the discretized model \DM{also is}. In addition and since the discretized model or Full Order Model (FOM) can be of large dimension, we recall the standard Loewner framework for the Model Order Reduction (MOR) using frequency domain interpolation. We recall the main steps to produce a Reduced Order Model (ROM) that approaches the FOM in a \DM{given} range of frequencies. We summarize the steps to follow in order to obtain a ROM that preserves the passive structure as well. Finally, we provide a constructive way to build a projector that allows to recover the physical meaning of the state variables from the ROM to the FOM. We use the one-dimensional wave equation and the Timoshenko beam as examples to show the versatility of the proposed approach. 
\end{abstract}


\begin{keyword}
Distributed port-Hamiltonian systems \sep Finite Element Method \sep Loewner framework \sep Structure-preserving discretization methods.



\end{keyword}

\end{frontmatter}


\section{INTRODUCTION}\label{sec:Intro}

The Hamiltonian formulation has been extended to distributed parameter systems in n-dimensional spaces with boundary energy flow in \citep{vanderSchaft2002JournalHamiltonian} and  proofs of existence and uniqueness of solutions have been presented for the one-dimensional case in \citep{LeGorrec2005JournalDirac}. This class of systems, called Boundary-Controlled Port-Hamiltonian Systems (BC-PHSs), has shown to be well suited for the modelling of open-loop infinite-dimensional systems with actuation and sensing located at the spatial boundaries, being applicable to the control of transport phenomena, waves, beams, and chemical reactors, among others \citep{Rashad2020JournalTwenty}.  

For numerical simulation of BC-PHSs, one has to discretize the infinite-dimensional part. Several techniques have been proposed using discrete exterior calculus \citep{Seslija2012JournalDiscrete}, mixed finite-elements \citep{Golo2004JournalHamiltonian}, finite volume \citep{Kotyczka2016ConferenceFinite}, finite-differences methods \citep{Trenchant2018JournalFinite}, among others. Recently, the Partitioned Finite Element
Method (PFEM) \citep{CardosoRibeiro2018ConferenceStructure}, \citep{Serhani2019ConferencePartitioned}, \citep{Cardoso2020JournalPartitioned} has shown a wide number of applications in one-, two- and three-dimensional domains \citep{Brugnoli2019JournalPort}, \citep{Haine2022ConferenceStructure}, allowing also to be applicable for the case of mixed boundary control \citep{Brugnoli2020ConferencePartitioned}, \citep{Brugnoli2022ConferenceExplicit}. In this paper, we particularize the PFEM to the  parametrized class of one-dimensional BC-PHSs introduced in \citep{LeGorrec2005JournalDirac}, covering the complete class of first-order spatial derivatives and all possible parametrizations of the boundary conditions. The matrices of the discretized model are explicitly defined in terms of the PDE matrices coefficients and the input/output parametrization matrices, being versatile to changes in the PDE structure and/or boundary conditions. Moreover, the proposed methodology guarantees that the discretized model preserves the initial energy structure of the PDE. That is, the Hamiltonian of the discretized model mimics the Hamiltonian of the initial PDE, showing that if the initial BC-PHS is {\it passive} (or {\it Impedance Energy Preserving} (IEP)), the discretized model is too.

The realization of the discretized model obtained using the PFEM contains only sparse matrices, which reduces the computational effort for time simulations. However, when the step of the spatial discretization is very small, the discretized model can contain a large number of internal variables, being expensive for numerical simulation and/or controller design. In this paper, we additionally recall the Loewner approach to find a Reduced Order Model (ROM) that captures the input/output behavior in a desired range of frequencies \citep{Mayo2007JournalFramework}, \citep{Antoulas2016ChapterATutorial}. This framework has shown to be an efficient tool for the reduction of infinite-dimensional systems \citep{Poussot2023ConferenceData}. Recently, it has been applied to distributed port-Hamiltonian systems \citep{Cherifi2021ConferenceApplication}  in 1D and \citep{Poussot2023ConferenceData} in 2D, and it is still a current research topic on the structure-preserving Model Order Reduction (MOR) field. In this paper, we combine the algorithms proposed in \citep{Benner2020JournalIdentification} \citep{Poussot2023ConferenceData},  to guarantee that the ROM remains passive. Moreover, we build a projector that allows us to reconstruct the original state of the discretized model from the ROM state. To the best of the authors knowledge this projection has not been exploited in the Loewner framework, so far. In this paper, we show that the projector can be straightforwardly constructed using the knowledge of the discretized model.

\subsection*{Summary of the contribution}
The PFEM is developped for the complete class of one-dimensional BC-PHSs with first-order spatial derivatives and all possible boundary inputs and outputs. The discretized model preserves the {\it passive} (or {\it impedance energy preserving}) structure. Additionally, a ROM technique is recalled using two recent algorithms in the Loewner framework. The ROM is guaranteed to be passive and a projector is obtained to reconstruct the state of the discretized model from the ROM state.

% complete methodology proposed in this paper can be summarized as follows: The class of BC-PHS is first discretized using the PFEM (first approximation). Then, the transfer function of the discretized model is interpolated in a desired frequency region using the Loewner approach  (second approximation). Finally, the interpolated model is reduced using a singular value decomposition (third and last approximation). For the second and third approximations, we recover a projection that allows to rebuild the state of the discretized model initially obtained using PFEM.

\subsection*{Paper oragnization}
The paper is organized as follows: in Section~\ref{Sec:PF}, we recall the structure of BC-PHSs and we emphasize the main contributions of this article. In Section~\ref{sec:PFEM}, the PFEM is presented. In Section~\ref{sec:Loewner}, the reduction of the discretized model is applied using the Loewner framework and in Section~\ref{Sec:Summary} a summary of the complete methodology is provided. The one-dimensional wave equation is used along the paper to exemplify the procedure, and in Section~\ref{Examples}, the Timoshenko beam is used to show the versatility of this approach. Finally, in Section~\ref{sec:Conclusion} some conclusions are drawn and future works are discussed.

%\newpage
\subsection*{Notation and Definitions}
Null square matrices of size $n \times n $ are denoted by $0_n$ and null rectangular matrices of size $n \times m$ by $0_{n \times m}$. For simplicity and when it is clear from the context, the subindex are omitted. The identity matrix of size $n$ is denoted by $I_n$. $u \in L_{loc}^2((0,\infty),\Real^n)$ means that the function $u(t)$ is locally square-integrable, {\it i.e.,} $\int_T \lvert u(t)\rvert _{\Real^n}^2 dt < +\infty $ for all compact subsets $T \subset \Real$.
\begin{definition}\label{Dissipativity}\citep[Definition~1.16.]{Villegas2007ThesisPort}
A state space system is called {\it dissipative} with respect to the supply rate $s(u(t),y(t))$ if there exists a function $S : X \rightarrow \mathbb{R}_+$, called the {\it storage function}, such that for all initial condition $x(0)=x_0 \in X$, all $t\geq 0$, and all input functions $u(t)$,
\begin{equation*}
\dot{S}(t) \leq s(u(t),y(t)).
\end{equation*} 
\end{definition}

\begin{definition}\label{Passivity}\citep[Definition~1.17.]{Villegas2007ThesisPort}
A state space system is called {\it passive} if it is {\it dissipative} with respect to the supply rate $s(u(t),y(t)) := u(t)^\top y(t)$ and is called {\it Impedance Energy Preserving} (IEP) in the particular case if $\dot{S}(t) := u(t)^\top y(t)$.
\end{definition}
\section{Problem Formulation}\label{Sec:PF}
\subsection{Boundary-controlled port-Hamiltonian systems in  a one-dimensional domain}
\label{sec:BC-PHS}
Boundary-controlled port-Hamiltonian systems are systems described by PDEs in which the state variables are chosen as the energy variables and the inputs and outputs are chosen as a linear combination of the co-energy variables evaluated at the boundaries of the spatial domain. In this section, we recall the main structure of this class of systems and the parametrization of the boundary conditions that lead to a well-posed BC-PHS.

\subsection{Partial Differential Equation}
We consider the following one-dimensional PDE:
\begin{equation}\label{Eq:PDE}
\begin{split}
\dfrac{\partial x}{\dt}(\zeta,t) &= P\,\dfrac{\partial e}{\partial \zeta} (\zeta,t) \,\DM{-}\, G\,e(\zeta,t), \quad  x(\zeta,0) = x_0(\zeta), \\
e(\zeta,t) &= \mathcal{H} (\zeta)\,x(\zeta,t),
 \end{split}
\end{equation}
%\DM{j'ai introduit le - devant le $G$, il faut maintenant le répercuter partout où $G$ apparaît dans le papier.}
with spatial domain $\zeta \in [a,b]$, time $t \geq0$, energy variable $x(\zeta,t) \in \mathbb{R}^n$, co-energy variable $e(\zeta,t)\in \mathbb{R}^n$,
%\begin{equation}\label{Co-Energy}
%e(\zeta,t) = \mathcal{H} (\zeta) x(\zeta,t),
%\end{equation} 
initial condition $x_0(\zeta)$, Hamiltonian density matrix $\mathcal{H}(\zeta)$, and matrices $P$, $G\in \Real^{n\times n}$. We consider the following assumtions:
\begin{itemize}
    \item $P$ is symmetric and invertible $P = P^\top$, 
    \item \DM{{$G+G^\top\geq 0$}, i.e. the symmetric part of $G$ is positive semi-definite},
    \item $\Ho(\zeta)$ is a bounded and continuously differentiable symmetric matrix-valued function of size $n \times n$ such that is positive definite for all $\zeta \in[a,b]$.
\end{itemize}
Without \DM{loss} of generality, we consider the following partitions of the matrices $P$, $G$ and $\mathcal{H}(\zeta)$:
\begin{equation*}\label{PDEMatrices}
\begin{split}
P &= \begin{pmatrix}
P_{11} & P_{12} \\
P_{21} & P_{22}
\end{pmatrix}, \,G = \begin{pmatrix}
G_{11} & G_{12} \\
G_{21}  & G_{22}
\end{pmatrix}, \, \Ho (\zeta) = \begin{pmatrix}
\Ho_1 (\zeta) & 0 \\
0 & \Ho_2 (\zeta)
\end{pmatrix}.
\end{split}
\end{equation*}
with $P_{11}$, $G_{11}$, $\Ho_1(\zeta) \in \Real^{n_1 \times n_1}$, $P_{12}$, $G_{12} \in \Real^{n_1 \times n_2}$, $P_{21}$, $G_{21} \in  \Real^{n_2 \times n_1}$ and $P_{22}$, $G_{22}$, $\Ho_2(\zeta) \in \Real^{n_2 \times n_2}$. Since $P$ is symmetric, $P_{11} = P_{11}^\top$, $P_{22} = P_{22}^\top$, and $P_{21} = P_{12}^\top$, and since $\Ho(\zeta)$ is symmetric positive definite  $\Ho_1(\zeta) = \Ho_1(\zeta)^\top >0$ and $\Ho_2(\zeta) = \Ho_2(\zeta)^\top >0$. The previous structure suggests the following energy and co-energy variables:
\begin{equation}
x(\zeta,t) =\begin{pmatrix}
x_1(\zeta,t) \\ x_2 (\zeta,t)
\end{pmatrix}, \quad e(\zeta,t) =\begin{pmatrix}
e_1(\zeta,t) \\ e_2 (\zeta,t)
\end{pmatrix}
\end{equation}
with $x_i(\zeta,t)$ and $e_i(\zeta,t) \in \Real^{n_i}$ for all $\zeta \in [a,b]$, $t\geq 0$, $i = \lbrace1,2\rbrace$, and $n_1+n_2 = n$.

\begin{remark}
When $n_1 = n_2$, the previous representation of the PDE \DM{reminds} the port-Hamiltonian formulation \DM{of mechanical systems}, where $x_1(\zeta,t)$ represents the set of general displacements and $x_2(\zeta,t)$ the set of generalized momenta. 
\end{remark}
\begin{remark}
\DM{
For the sake of simplicity, in this paper the standard notations, originally introduced in~\citep{LeGorrec2005JournalDirac} have been adapted in 
the PDE \eqref{Eq:PDE}:  $P$ stands for $P_1$,  and  $G$ now stands for $G_0-P_0$; thus {$G+G^\top\geq 0$} is an equivalent formulation for $G_0=G_0^\top \geq 0$.
}
\end{remark}

\subsection{Hamiltonian}

The Hamiltonian associated to the PDE in \eqref{Eq:PDE} is defined by
\begin{equation}\label{Eq:Hamiltonian}
H(t) := \dfrac{1}{2} \int_a^b {x(\zeta,t)^\top \mathcal{H}(\zeta) x(\zeta,t)}d\zeta = \dfrac{1}{2} \int_a^b {x^\top e}d\zeta. 
\end{equation}
Using the fact that $P = P^\top$ and $G+G^\top\geq 0$, the energy balance is obtained as follows:
\begin{equation}
\begin{split}
\dot{H}(t) &= \dfrac{1}{2} \int_a^b \left( \dot{x}^\top e + x^\top \dot{e} \right) d\zeta, \\
&= \dfrac{1}{2} \int_a^b \left( (P \tfrac{\partial e}{\partial \zeta}-Ge)^\top e + x^\top\Ho (P \tfrac{\partial e}{\partial \zeta}-Ge) \right) d\zeta, \\
&= \dfrac{1}{2} \int_a^b \left( (P \tfrac{\partial e}{\partial \zeta}-Ge)^\top e + e^\top (P \tfrac{\partial e}{\partial \zeta}-Ge) \right) d\zeta, \\
&\leq   \dfrac{1}{2} \left[e^\top P e \right]_a^b.
\end{split} 
\end{equation}
In the following, we define the inputs and outputs for the PDE \eqref{Eq:PDE}, such that the Hamiltonian \eqref{Eq:Hamiltonian} qualifies as a storage function for {\it passivity} as defined in Definition \ref{Passivity}.

\subsection{Boundary inputs and outputs}\label{SubSec:BIO}
In what follows we define boundary inputs and outputs, and derive the associated energy balance.
\begin{theo}%[\cite{LeGorrec2005JournalDirac}, Theorem 4.2 and Theorem 4.4]
\label{Theo:Yann}
Let %$\Sigma = \left( \begin{smallmatrix}0_n & I_n \\ I_n & 0_n \end{smallmatrix} \right)$ and
$V_\mathcal{B},V_\mathcal{C} \in \Real^{n \times 2n}$ be two full-rank matrices such that 
\begin{equation}
V_i \left(\begin{matrix}P^{-1}&0\\0&-P^{-1}\end{matrix}\right) V_i^\top = 0_n, \quad i=\lbrace \mathcal{B}, \mathcal{C} \rbrace.
\end{equation}
Define the {\it boundary input} and {\it boundary output} as
\begin{equation}\label{InputOutput}
u(t) = V_\mathcal{B}\begin{pmatrix}
e(b,t) \\
e(a,t)
\end{pmatrix}, \quad y(t) = V_\mathcal{C}\begin{pmatrix}
e(b,t) \\
e(a,t)
\end{pmatrix}.
\end{equation}
The system \eqref{Eq:PDE}-\eqref{InputOutput} is a boundary control system. Furthermore, for all inputs $u(t) \in L_{loc}^2((0,\infty),\Real^n)$, initial \DM{data}  $x_0(\zeta)\in H^{1}((a,b);\Real^n)$ with $u(0) = V_\mathcal{B}\left(\begin{smallmatrix}
e(b,0) \\
e(a,0)
\end{smallmatrix}\right)$, the following balance equation with respect to the Hamiltonian \eqref{Eq:Hamiltonian} is satisfied:
\begin{equation}\label{Balance}
\dot{H}(t) \leq u(t) ^\top y(t).
\end{equation}
\end{theo}

\begin{Proof}
The proof follows directly from
\citep[Theorem 4.4]{LeGorrec2005JournalDirac}, considering $W=V_\mathcal{B} \mathcal{R}^{-1}$ and $\tilde{W}=V_\mathcal{C} \mathcal{R}^{-1}$ with $\mathcal{R}:=\frac{1}{\sqrt{2}}\left(\begin{smallmatrix}P & -P \\ I & I\end{smallmatrix} \right)$
\end{Proof}

%Notice that existence and uniqueness of solution of the system \eqref{Eq:PDE}-\eqref{InputOutput} is shown for the case in which $u(t) = 0$ \citep[Theorem 4.1]{LeGorrec2005JournalDirac}. 
%Moreover, if {$G+G^\top = 0$}, the differential operator is the generator of a unitary semigroup on $L^2((a,b); \Real^n)$ and the balance becomes $\dot{H}(t) = u(t) ^\top y(t)$, implying an IEP structure. In the following lemma, we recall some properties of the input and output matrices. 

Note that it has been shown in \citep[Theorem 4.1]{LeGorrec2005JournalDirac} that if {$G+G^\top = 0$}, the differential operator in \eqref{Eq:PDE}-\eqref{InputOutput} is the generator of a unitary semigroup on $L^2((a,b); \Real^n)$ and the balance equation on the energy reads $\dot{H}(t) = u(t) ^\top y(t)$, implying an IEP structure. In what follows we recall some properties of the input and output matrices. 

\begin{lemma}\label{Lemma:InputOutput}
The matrices $V_\mathcal{B}$ and $V_\mathcal{C}$ introduced in Theorem \ref{Theo:Yann} satisfy
\begin{equation}
V_\mathcal{C}^\top V_\mathcal{B} = \blue{\dfrac{1}{2}} \left(\begin{matrix}P&0 \\ 0 & -P\end{matrix}\right) - \tilde{J}_\xi
\end{equation}
for some $2n \times 2n$ skew-symmetric matrix $\tilde{J}_\xi = -\tilde{J}_\xi^\top$.
\end{lemma}
\begin{Proof}
The proof follows directly from \citep[Lemma~A.1]{Macchelli2020JournalExponential} replacing $W$ by $V_\mathcal{B}\mathcal{R}^{-1}$ and $\tilde{W}$ by $V_\mathcal{C}\mathcal{R}^{-1}$ where $\mathcal{R}$ is defined in the proof of Theorem \ref{Theo:Yann}.
\end{Proof}

\begin{example}\label{Example}
(One-dimensional wave equation) Let us consider the following 1D wave equation: 
\begin{equation*}
\rho (\zeta)\dfrac{\partial^2 w}{\partial t^2} (\zeta,t)  = \dfrac{\partial }{\partial \zeta} \left( T(\zeta)\dfrac{\partial w}{\partial \zeta }(\zeta,t)\right), \;
w(\zeta,0) = w_0(\zeta), \;
\dfrac{\partial w}{\partial t}(\zeta,0) = v_0(\zeta),
\end{equation*}
where $\zeta \in [a,b]$ is the space, $t \geq 0$ is the time, $w(\zeta,t) \in \mathbb{R}$ is the wave deformation, $w_0(\zeta) \in H^1(a,b)$ is the initial deformation and $v_0(\zeta) \in L^2(a,b) $ is the initial velocity. $\rho(\zeta)$ and $T(\zeta)$ are both positive definite functions for all $\zeta \in [a,b]$ representing the mass density and the modulus of elasticity, respectively. Defining $x_1(\zeta,t) := \dfrac{\partial w}{\partial \zeta} (\zeta,t)$, $x_2(\zeta,t) := \rho(\zeta)\dfrac{\partial w}{\partial t} (\zeta,t)$, $e_1(\zeta,t) := T(\zeta)x_1(\zeta,t)$, and $e_2(\zeta,t): = \dfrac{1}{\rho(\zeta)}x_2(\zeta,t)$, the previous PDE can be written as in~\eqref{Eq:PDE} as follows:
\begin{equation*}\label{Eq:1DWavePHS}
\dfrac{\partial }{\partial t}\begin{pmatrix}
x_1 \\ x_2
\end{pmatrix} = \begin{pmatrix}
0 & 1 \\ 1 & 0
\end{pmatrix}\dfrac{\partial}{\partial \zeta}\begin{pmatrix}
e_1 \\ e_2
\end{pmatrix},\; 
x_1(\zeta,0) = \dfrac{dw_0}{d\zeta}(\zeta), \;
 x_2(\zeta,0) = \rho(\zeta)v_0(\zeta), \\
\end{equation*} 
with ${P}_{12} ={P}_{21}= 1$, ${P}_{11} ={P}_{22}= {G}_{12} = {G}_{21}={G}_{11} ={G}_{22}= 0$, $\mathcal{H}_1(\zeta) = T(\zeta)$ and $\mathcal{H}_2(\zeta) = \rho(\zeta)^{-1}$. Different parameterizations of the inputs and outputs can be formulated from the previous representation. This is the case of Neumann boundary inputs and velocity outputs
\begin{equation}\label{Neumann}
V_\mathcal{B} = \begin{pmatrix}
0 & 0 & 1 & 0 \\ 
1 & 0 & 0 & 0 
\end{pmatrix}, \quad V_\mathcal{C} = \begin{pmatrix}
0 & 0 & 0 & -1 \\ 
0 & 1 & 0 & 0 
\end{pmatrix},
\end{equation}
or mixed boundary inputs and outputs
\begin{equation}\label{mixed}
V_\mathcal{B} = \begin{pmatrix}
0 & 0 & 0 & 1 \\ 
1 & 0 & 0 & 0 
\end{pmatrix}, \quad V_\mathcal{C} = \begin{pmatrix}
0 & 0 & -1 & 0 \\ 
0 & 1 & 0 & 0 
\end{pmatrix},
\end{equation}
for instance. In both cases, the balance equation \eqref{Balance} is satisfied. 
\end{example}

\begin{remark}
A more general formulation of the BC-PHS \eqref{Eq:PDE}-\eqref{InputOutput} is presented in \citep{LeGorrec2005JournalDirac}, \DM{which includes higher-order differential operators, thus  allowing to tackle such cases as the Euler-Bernoulli beam}. For simplicity and clarity of the explanation of this article, we only deal with the case of first-order differential operators. 
\iffalse
\DM{Note that the application of PFEM to the Euler-Bernoulli beam has been presented in detail in \citep[\S 6.3, pp. 527--530]{Cardoso2020JournalPartitioned}; one of the features of 1D higher-order
operators is that many more variables are involved at the boundary; typically for the Euler-Bernoulli beam,
we end up with 4 boundary components instead of only 2 for the vibrating string: boundary values as well
as spatial derivatives of these boundary values have to be taken into account. Moreover on the numerical side, for the choice of the Finite Element families, since a higher regularity in space is required, we resort to Hermite polynomials instead of Lagrange polynomials.
}
\fi
\end{remark}

\subsection{Proposed methodology}%Main contributions of this article}
\begin{figure}[!h]
\begin{center}
\includegraphics[width=0.4\textwidth]{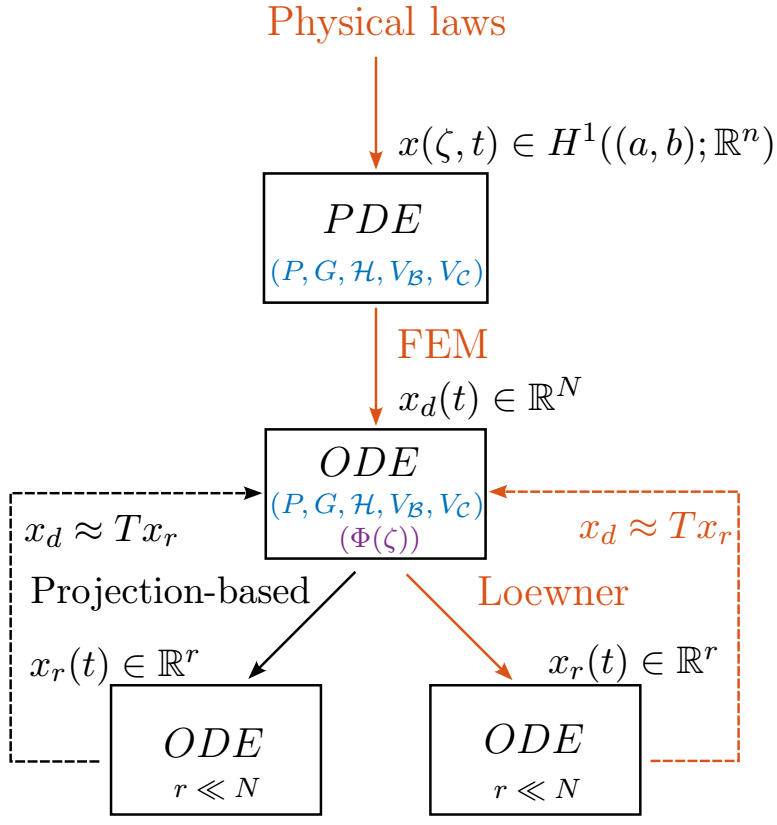}  
\caption{Diagram of the article contributions. Read as follows: from the physical laws, a class of boundary control systems is modeled using PDEs. The PDEs have state $x(\zeta,t)$ with $\zeta \in[a,b]$ and $t\geq 0$. The dynamics of the PDE is completely described by the matrices \blue{$(P,G,\Ho,V_\mathcal{B},V_\mathcal{C})$}. The spatial variable $\zeta$ is discretized using the FEM obtaining an $ODE$ with state $x_d(t)$ of size $N$. The matrices of the obtained ODE are defined as a function of the PDE matrices \blue{$(P,G,\Ho,V_\mathcal{B},V_\mathcal{C})$} and the basis functions \violet{$\Phi(\zeta)$} (required for  discretization). Finally, the ODE is reduced using the Loewner framework, obtaining another $ODE$ with state $x_r(t)$ with size $r \ll N$. We provide the projector $'T'$ that enables to approximate the state $x_d(t)$ of the FOM using the state $x_r(t)$ of the ROM.   } 
\label{Fig:Diagram}                           
\end{center}
\end{figure}

In Fig \ref{Fig:Diagram}, we show a diagram of the methodology proposed in this paper. The main steps are described hereafter : %{\bf main contributions} are summarized as follows:
\begin{itemize}
\item The BC-PHS \eqref{Eq:PDE}-\eqref{InputOutput} is discretized in such a way that the dynamics of the discretized model is defined as a function of the PDE matrices \blue{$(P,G,\Ho,V_\mathcal{B},V_\mathcal{C})$}. This means that the obtained ODE can be tuned for different types of PDEs and different types of boundary conditions. 
\item The passive structure of the BC-PHS \eqref{Eq:PDE}-\eqref{InputOutput} is preserved at the approximation level \DM{(on the FOM given by the ODE)}, {\it i.e.,} the balance equation \eqref{Balance} is exactly preserved, meaning that if \eqref{Eq:PDE}-\eqref{InputOutput} is {\it passive} (or {\it IEP}), the discretized model remains {\it passive} (or {\it IEP}). 
\item The Loewner framework is recalled and summarized in such a way that we are able to obtain a passive ROM. Moreover, in order to preserve the physical meaning of the variables, in this paper, similar to the projection-based MOR techniques, we provide a projector $'T'$ that allows to approximate the state of the FOM from the state of the ROM.   
\end{itemize}

\section{Partitioned Finite Element Method}\label{sec:PFEM}

In this section, we project the PDE \eqref{Eq:PDE}-\eqref{InputOutput} \DM{onto} a finite-dimensional space. Then, we show that if the BC-PHS \eqref{Eq:PDE}-\eqref{InputOutput} is {\it passive} (or {\it IEP}), the finite-dimensional model preserves this structure with the storage energy given by the projected Hamiltonian derived from \eqref{Eq:Hamiltonian}. For simplicity of the explanation, in this section we develop the projected version of the weak formulation, which is developped in \ref{Appendix_WF}. 

\subsection{Spatial discretization through the projections.}\label{Sec:Projections}

The PDE \eqref{Eq:PDE} is written in the following form:
\begin{equation}\label{PDE_and_Constitutive}
\begin{cases}
\dot{x}_1 = P_{11}\partial_\zeta e_1 + P_{12}\partial_\zeta e_2 - G_{11} e_1 - G_{12} e_2, \\
\dot{x}_2 = P_{21}\partial_\zeta e_1 + P_{22}\partial_\zeta e_2 - G_{21} e_1 - G_{22} e_2, \\
e_1 = \Ho_1 x_1, \\
e_2 = \Ho_2 x_2,
\end{cases}
\end{equation}
in which we have used the notation $\dot{f}$ for $\tfrac{\partial f}{\partial t}$ and $\partial_\zeta f$ for $\tfrac{\partial f}{\partial \zeta}$. We define the scalar spatial-dependent functions $\phi_1^{i_1}(\zeta) \in V_1 \subset H^1(a,b)$ and $\phi_2^{i_2}(\zeta) \in V_2 \subset H^1(a,b)$ for $i_1 = \lbrace 1,\cdots , N_1 \rbrace$ and $i_2 = \lbrace 1,\cdots , N_2 \rbrace$, respectively. The vector functions are defined as:
\[\phi_1(\zeta)= \begin{pmatrix}
\phi_1^1(\zeta) \\ \vdots \\ \phi_1^{N_1}(\zeta)
\end{pmatrix}, \quad 
\phi_2(\zeta)= \begin{pmatrix}
\phi_2^1(\zeta) \\ \vdots \\ \phi_2^{N_2}(\zeta)
\end{pmatrix}.\]
 Then, all state variables contained in $x_1$ and $e_1$ are projected through $\phi_1(\zeta)$, whereas the ones contained in $x_2$ and $e_2$ through $\phi_2(\zeta)$, {\it i.e.,}
 \begin{equation}\label{ApproxVariables}
 \begin{split}
 x_i(\zeta,t)& \approx x_i^{ap}(\zeta,t) := \Phi_i(\zeta)^\top x_{id}(t),\\
  e_i(\zeta,t)& \approx e_i^{ap}(\zeta,t) := \Phi_i(\zeta)^\top e_{id}(t),
%   x_2(\zeta,t)& \approx x_2^{ap}(\zeta,t) := \Phi_2(\zeta)^\top x_{2d}(t),\\
%  e_2(\zeta,t)& \approx e_2^{ap}(\zeta,t) := \Phi_2(\zeta)^\top e_{2d}(t),
  \end{split}
 \end{equation}
 with $i = \lbrace{1,2\rbrace}$, $x_i^{ap}$, $e_i^{ap} \in \Real^m$ being the approximated vector functions, 
\[
\Phi_1= \begin{pmatrix}
\phi_1  & \cdots & 0 \\ 
 \vdots & \ddots & \vdots \\ 
0 &  \cdots & \phi_1 \\ 
\end{pmatrix}_{n_1N_1 \times n_1} , 
\quad \Phi_2 = \begin{pmatrix}
\phi_2  & \cdots & 0 \\ 
 \vdots & \ddots & \vdots \\ 
0 &  \cdots & \phi_2 \\ 
\end{pmatrix}_{n_2N_2 \times n_2} ,
\]
being the projections matrices, and $x_{1d}$, $e_{1d} \in \Real^{n_1N_1}$ and $x_{2d}$, $e_{2d} \in \Real^{n_2N_2}$ the time domain variables. We have omitted the spatial and temporal dependence since it is clear from the context.

The objective is to discretize the spatial variable $\zeta$ of all equations in \eqref{PDE_and_Constitutive}. To this end, from \eqref{PDE_and_Constitutive}, we multiply the first and third equations by $\Phi_1$ on the left and the second and forth equations by $\Phi_2$ on the left. Then, we integrate over the spatial domain $\zeta \in [a,b]$ the four equations. Then, \eqref{PDE_and_Constitutive} becomes:

\begin{equation}\label{Eq:ODE1}
\begin{cases}
E_1 \dot{x}_{1d} = (D_{11}^P-D_{11}^G)e_{1d} + (D_{12}^P-D_{12}^G)e_{2d}, \\
E_2 \dot{x}_{2d} = (D_{21}^P-D_{21}^G)e_{1d} + (D_{22}^P-D_{22}^G)e_{2d}, \\
E_1 e_{1d} = Q_1 x_{1d}, \\
E_2 e_{2d} = Q_2 x_{2d}, \\
\end{cases}
\end{equation}
with square matrices
\begin{align*}
E_i = \int_a^b \Phi_i \Phi_i ^\top d\zeta, \quad 
Q_i = \int_a^b \Phi_i \Ho_i \Phi_i ^\top d\zeta, \quad i = \lbrace 1,2 \rbrace,
\end{align*}
and for $i = \lbrace 1,2 \rbrace$ and $j = \lbrace 1,2 \rbrace$ rectangular matrices
\begin{equation*}
D_{ij}^P = \int_a^b \Phi_i {P}_{ij} (\partial_\zeta \Phi_j) ^\top d\zeta, \quad  D_{ij}^G = \int_a^b \Phi_i {G}_{ij} \Phi_j ^\top d\zeta. \\
\end{equation*}
\begin{comment}
\begin{equation*}
\begin{array}{lccl}
D_{11}^P = \int_a^b \Phi_1 {P}_{11} (\partial_\zeta \Phi_1) ^\top d\zeta, & & & D_{11}^G = \int_a^b \Phi_1 {G}_{11} \Phi_1 ^\top d\zeta,\\
D_{12}^P = \int_a^b \Phi_1 {P}_{12} (\partial_\zeta \Phi_2) ^\top d\zeta, & & & D_{12}^G = \int_a^b \Phi_1 {G}_{12}  \Phi_2 ^\top d\zeta,\\
D_{21}^P = \int_a^b \Phi_2 {P}_{21} (\partial_\zeta \Phi_1) ^\top d\zeta, & & & D_{21}^G = \int_a^b \Phi_2 {G}_{21} \Phi_1 ^\top d\zeta,\\
D_{22}^P = \int_a^b \Phi_2 {P}_{22} (\partial_\zeta \Phi_2) ^\top d\zeta, & & & D_{22}^G = \int_a^b \Phi_2 {G}_{22}  \Phi_2 ^\top d\zeta.
\end{array}
\end{equation*}
\end{comment}
The equations in \eqref{Eq:ODE1} can be compactly written as written as:
\begin{equation}\label{Eq:ODE2}
\begin{cases}
\begin{pmatrix}
E_1 & 0 \\ 0 & E_2
\end{pmatrix} \begin{pmatrix}
 \dot{x}_{1d} \\  \dot{x}_{2d}
\end{pmatrix} = \begin{pmatrix}
D_{11}^P-D_{11}^G & D_{12}^P-D_{12}^G \\
D_{21}^P-D_{21}^G & D_{22}^P-D_{22}^G
\end{pmatrix} \begin{pmatrix}
e_{1d} \\ e_{2d}
\end{pmatrix}, \\
\begin{pmatrix}
E_1 & 0 \\ 0 & E_2
\end{pmatrix} \begin{pmatrix}
{e}_{1d} \\  {e}_{2d}
\end{pmatrix} = \begin{pmatrix}
Q_1 & 0 \\
0 & Q_2
\end{pmatrix} \begin{pmatrix}
x_{1d} \\ x_{2d}
\end{pmatrix},
\end{cases}
\end{equation}
Since the system under study \eqref{Eq:PDE}-\eqref{InputOutput} is an infinite-dimensional port-Hamiltonian system, we aim to find a finite-dimensional port-Hamiltonian representation for the equations in \eqref{Eq:ODE2} with conjugated inputs and outputs, {\it i.e.,} we aim to find the following representation:
\begin{equation}\label{StateSpace2}
\Sigma_d\begin{cases}
E \dot{x}_d = (J-R) e_d + B u , \quad x_d(0) = x_{d0}, \\
E e_d = Qx_d, \\
\quad  y = B^\top e_d.
\end{cases}
\end{equation}
with $x_d = (x_{1d}^\top , x_{2d}^\top)^\top$, $e_d = (e_{1d}^\top , e_{2d}^\top)^\top$, $E$ a mass matrix (that is symmetric and positive-definite), $Q$ an energy matrix, $J$ skew-symmetric, $R$ symmetric and positive semi-definite, and $B$ the input matrix. We can notice that in the state-space representation \eqref{Eq:ODE2}, the input $u(t)$ and output $y(t)$ do not appear explicitly. Then, we define the input matrix as the power conjugated of the output matrix defined in \eqref{InputOutput}. To this end, we project the input and output as follows:
\begin{equation}
\begin{split}
u(t)&= V_\mathcal{B} \begin{pmatrix}
e(b,t) \\
e(a,t)
\end{pmatrix} \approx V_\mathcal{B} \begin{pmatrix}
e_1^{ap}(b,t) \\
e_2^{ap}(b,t) \\
e_1^{ap}(a,t) \\
e_2^{ap}(a,t)
\end{pmatrix} = V_\mathcal{B} \Omega _{ba}^\top \begin{pmatrix}
e_{1d} (t) \\
e_{2d}(t)
\end{pmatrix}, \\
y(t)&= V_\mathcal{C} \begin{pmatrix}
e(b,t) \\
e(a,t)
\end{pmatrix} \approx V_\mathcal{C} \begin{pmatrix}
e_1^{ap}(b,t) \\
e_2^{ap}(b,t) \\
e_1^{ap}(a,t) \\
e_2^{ap}(a,t)
\end{pmatrix} = V_\mathcal{C} \Omega _{ba}^\top \begin{pmatrix}
e_{1d} (t) \\
e_{2d}(t)
\end{pmatrix}, 
\end{split}
\end{equation}
with 
\begin{equation}
    \Omega _ {ba} := \begin{pmatrix}
\Phi_1(b) & 0_{n_1N_1 \times n_2} & \Phi_1(a) & 0_{n_1N_1 \times n_2} \\
0_{n_2N_2 \times n_1} & \Phi_2(b) & 0_{n_2N_2 \times n_1} & \Phi_2(a)  \\
\end{pmatrix}
\end{equation}
and then, we define the input matrix as conjugated with respect to the output one
\begin{equation}
B = \Omega_{ba}V_\mathcal{C}^\top.
\end{equation}
Then, by adding and subtracting $\pm Bu(t)$ at the right-hand side equation of the first block equation in \eqref{Eq:ODE2} and using the input projection $u \approx V_\mathcal{B}\Omega_{ba}^\top e_d$, we obtain the state space representation \eqref{StateSpace2} with
%\begin{equation}\label{StateSpace2}
%\begin{split}
%E \dot{x}_d (t) &= F e_d(t) + Bu(t), \quad x_d(0) = x_{d0}, \\
%y(t) &= B^\top e_d(t). \end{split}
%\end{equation}
\begin{equation}\label{Eq:JE}
\begin{split}
E &= \begin{pmatrix} {E}_1 & 0_{n_1 N_1 \times n_2 N_2} \\ 0_{n_2 N_2 \times n_1 N_1} & {E}_2 \end{pmatrix}, \\ Q &= \begin{pmatrix} {Q}_1 & 0_{n_1 N_1 \times n_2 N_2} \\ 0_{n_2 N_2 \times n_1 N_1} & {Q}_2 \end{pmatrix}, \\
J &= D^P - \dfrac{1}{2}(D^G-{D^G}^\top)- \Omega_{ba} V_\mathcal{C}^\top V_\mathcal{B} \Omega_{ba}^\top, \\
R &=  \dfrac{1}{2}(D^G+{D^G}^\top), 
\end{split}
\end{equation}
\begin{equation*}
D^P = \begin{pmatrix}
D_{11}^P & D_{12}^P \\
D_{21}^P & D_{22}^P
\end{pmatrix}, \quad D^G = \begin{pmatrix}
D_{11}^G & D_{12}^G \\
D_{21}^G & D_{22}^G
\end{pmatrix}
\end{equation*}
and initial data
\[ x_{d0} = \begin{pmatrix}
x_{1d}(0)\\
x_{2d}(0)
\end{pmatrix} = \begin{pmatrix}
    E_1^{-1} \int_a^b \Phi_1 x_1(\zeta,0)d\zeta\\
    E_2^{-1} \int_a^b \Phi_2 x_2(\zeta,0)d\zeta
\end{pmatrix}.\] %\quad \begin{pmatrix} x_1(\zeta,0) \\ x_2(\zeta,0)\end{pmatrix}= \begin{pmatrix}\Phi(\zeta) x_{1d}(0) \\ \Phi(\zeta) x_{2d}(0)  \end{pmatrix}.\]
In order to preserve the port-Hamiltonian structure of the state-space representation \eqref{StateSpace2}, we have to guarantee the skew-symmetry of the matrix $J$ and the positive semi-definiteness of $R$. Notice that the matrix $J$ depends on $P$, $G$, $V_\mathcal{B}$ and $V_\mathcal{C}$ of the PDE \eqref{Eq:PDE}, whereas $R$ depends on $G$ only. We use the following properties for the main result:
\begin{proposition}\label{MatrixProperties}
The following properties are satisfied:
\begin{equation}\label{Eq:D1Prop}
\begin{split}
(a) \quad D_{12}^P + {D_{21}^P}^\top &= \left[\Phi_1 P_{12}\Phi_2^\top\right]_a^b, \\
(b) \quad  D_{11}^P + {D_{11}^P}^\top &= \left[\Phi_1 P_{11}\Phi_1^\top\right]_a^b, \\
(c) \quad  D_{22}^P + {D_{22}^P}^\top &= \left[\Phi_2 P_{22}\Phi_2^\top\right]_a^b.
\end{split}
\end{equation}
\end{proposition}
\begin{Proof}
Property $(a)$ writes:
\begin{equation*}
\begin{split}
 D_{12}^P + {D_{21}^P}^\top &= \int_a^b \Phi_1 {P}_{12} (\partial_\zeta \Phi_2) ^\top d\zeta +  \int_a^b (\partial_\zeta \Phi_1) {P_{21}}^\top \Phi_2 ^\top d\zeta, \\
  & = \int_a^b \frac{\partial }{\partial \zeta} \left( \Phi_1 P_{12} \Phi_2 ^\top \right)d\zeta = \left[ \Phi_1 {P_{12}}  \Phi_2^\top \right]_a^b,\\
\end{split}
\end{equation*}
where we have used the fact that $P_{21} = P_{12}^\top$. The same line follows for $(b)$ and $(c)$ using the facts that $P_{11} = P_{11}^\top$ and $P_{22} = P_{22}^\top$, respectively.
\end{Proof}

\begin{proposition}\label{Prop:JR}
The matrix $J$ defined in \eqref{Eq:JE} is skew-symmetric, and  matrix $R$ is positive semi-definite.
\end{proposition}
\begin{Proof}
The term $-\tfrac{1}{2}(D^G-{D^G}^\top)$ from $J$ in \eqref{Eq:JE} is clearly skew-symmetric. It remains to show that the remaining part is also skew-symmetric, {\it i.e.,} $J_{P\Bo \Co} := D^P- \Omega_{ba} V_\mathcal{C}^\top V_\mathcal{B} \Omega_{ba}^\top$ is such that $J_{P\Bo \Co} + J_{P\Bo \Co}^\top = 0$. Notice that, using Proposition \ref{MatrixProperties}, one can obtain
\[
\begin{split}
D^P + {D^P}^\top &= \begin{pmatrix}
[\Phi_1 P_{11}\Phi_1^\top]_a^b & [\Phi_1 P_{12}\Phi_2^\top]_a^b \\
[\Phi_2 P_{21}\Phi_1^\top]_a^b & [\Phi_2 P_{22}\Phi_2^\top]_a^b
\end{pmatrix}= \Omega_{ba} \begin{pmatrix}
    P & 0 \\
    0 & -P
\end{pmatrix} \Omega_{ba}^\top.
\end{split}
\]
Then 
\[J_{P\Bo \Co} + J_{P\Bo \Co}^\top = \Omega_{ba}\left( \begin{pmatrix}
    P & 0 \\ 0 & -P
\end{pmatrix} - V_\Co ^\top V_\Bo - V_\Bo ^\top V_\Co \right) \Omega_{ba}^\top
\]
Finally, we replace $V_\Co ^\top V_\Bo$ and $V_\Bo ^\top V_\Co$ using Lemma \ref{Lemma:InputOutput} to obtain
\[J_{P\Bo \Co} + J_{P\Bo \Co}^\top = \Omega_{ba}\left( \tilde{J}_\xi -  \tilde{J}_\xi \right) \Omega_{ba}^\top = 0.
\]
Now, we show that $R$ is positive semi definite. To this end, we compute every block element of the matrix $R$ as follows:
\[
\begin{split}
    D_{11}^G + {D_{11}^G} ^\top & = \int_a ^b \Phi_1 \left( G_{11} + G_{11}^\top \right)\Phi_1 ^\top d\zeta, \\
    D_{22}^G + {D_{22}^G} ^\top & = \int_a ^b \Phi_2 \left( G_{22} + G_{22}^\top \right)\Phi_2 ^\top d\zeta, \\
    D_{12}^G + {D_{21}^G} ^\top & = \int_a ^b \Phi_1 \left( G_{12} + G_{21}^\top \right)\Phi_2 ^\top d\zeta, \\
    D_{21}^G + {D_{12}^G} ^\top & = \int_a ^b \Phi_2 \left( G_{21} + G_{12}^\top \right)\Phi_1 ^\top d\zeta, \\
\end{split}
\]
Then 
\[
\begin{split}
R &= \frac{1}{2}\int_a ^b \begin{pmatrix}
    \Phi_1 & 0 \\ 0 & \Phi_2
\end{pmatrix}\begin{pmatrix}
    G_{11}+G_{11}^\top & G_{12}+G_{21}^\top \\ G_{21}+G_{12}^\top & G_{22}+G_{22}^\top
\end{pmatrix}\begin{pmatrix}
    \Phi_1^\top & 0 \\ 0 & \Phi_2^\top
\end{pmatrix}d\zeta, \\
&= \frac{1}{2}\int_a ^b \begin{pmatrix}
    \Phi_1 & 0 \\ 0 & \Phi_2
\end{pmatrix} \left(G + G^\top \right)\begin{pmatrix}
    \Phi_1^\top & 0 \\ 0 & \Phi_2^\top
\end{pmatrix} d\zeta. \\
\end{split}
\]
Since by definition $G+G^\top \geq 0$, hence $R\geq 0$.
\end{Proof}

\subsection{Structure Preservation of the Hamiltonian}
We define the Hamiltonian of the discretized model replacing the approximations \eqref{ApproxVariables} in the Hamiltonian function \eqref{Eq:Hamiltonian}. The discretized Hamiltonian reads:
\begin{equation}\label{DiscretizedHamiltonian}
H_d(t) = \frac{1}{2} x_d(t) ^\top E e_d(t) = \frac{1}{2} x_d(t) ^\top Q x_d(t).
\end{equation}

\begin{proposition}
According to Definition \ref{Passivity}, the discretized model in \eqref{StateSpace2} is {\it passive} with respect to the storage energy function \eqref{DiscretizedHamiltonian}. Moreover, if the BC-PHS \eqref{Eq:PDE}-\eqref{InputOutput} is IEP, the discretized model in \eqref{StateSpace2} remains IEP.  
\end{proposition}
\begin{Proof}
The time derivative of the discretized Hamiltonian $H_d$ in \eqref{DiscretizedHamiltonian} writes:
\begin{equation}\label{Eq:dHddt}
\begin{split}
\dot{H}_d(t) &= x_d^\top Q \dot{x}_d, \\
&= x_d^\top QE^{-1}  E\dot{x}_d, \\
&= e_d^\top EE^{-1}  \left( (J-R)e_d + Bu \right), \\
&= e_d^\top  \left( (J-R)e_d + Bu \right), \\
&= e_d^\top    (J-R)e_d + e_d^\top  Bu , \\
&= -e_d^\top R e_d + e_d^\top Bu , \\
&= -e_d^\top R e_d +y^\top u.
\end{split} 
\end{equation}
Since $R\geq 0$ (See Proposition \ref{Prop:JR}), $\dot{H}_d\leq y^\top u$, hence $\Sigma_d$ in \eqref{StateSpace2} is {\it passive}. Now, we show that if \eqref{Eq:PDE}-\eqref{InputOutput} is IEP, the discretized model \eqref{StateSpace2} is IEP. As pointed out in subsection \ref{SubSec:BIO}, if $G+G^\top = 0$, system \eqref{Eq:PDE}-\eqref{InputOutput} is IEP. Notice that $R$ defined in \eqref{Eq:JE} is dependent on the $G$ matrix, only. Moreover, as it is developed in the proof of Proposition \ref{Prop:JR}, the matrix $R$ is a null matrix if $G+G^\top$ is null, {\it i.e.,} if $G+G^\top = 0$, then $R = 0$. Then, if \eqref{Eq:PDE}-\eqref{InputOutput} is IEP $(G+G^\top = 0)$, then $R = 0$ and \eqref{Eq:dHddt} becomes $\dot{H}_d = y ^\top u$. This concludes the proof.
\end{Proof}

\begin{example}\label{Example2}(One-dimensional wave equation)
We refer to Example \ref{Example} for the PDE formulation. For simplicity of the explanation, we consider $\phi_1(\zeta) =\phi_2(\zeta)=:\phi(\zeta) $. As shown in Example \ref{Example}, the PDE is partitioned in two states of size one, {\it i.e.,} $n_1=n_2 = 1$. Then, $\Phi_1(\zeta) = \phi_1(\zeta)= \phi(\zeta)$ and $\Phi_2(\zeta) = \phi_2(\zeta) = \phi(\zeta)$. In Figure \ref{Fig:HatFunctions}, we show one example of the selected basis functions, which are stacked in the vector $\phi(\zeta)$ with $N_1 = N_2 = N = 4$.
\begin{figure}[!h]
\begin{center}
\includegraphics[width=0.48\textwidth]{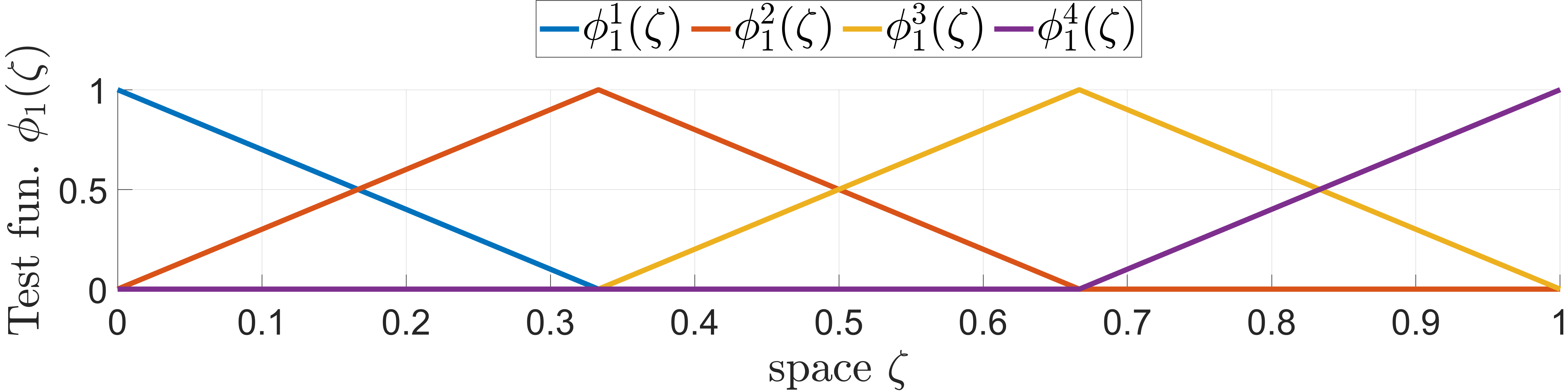}  
\caption{Example of basis functions using $N=4$.} 
\label{Fig:HatFunctions}     
\end{center}
\end{figure}
For simplicity, we consider constant densities $T(\zeta) = T_0$ and $\rho(\zeta) = \rho_0$. Then, the matrices of the system \eqref{StateSpace2}-\eqref{Eq:JE} are as follows:
\begin{equation}
    E_1 =E_2 = \tfrac{h}{6} \left(\begin{smallmatrix}
        2 & 1 & 0 & 0\\
        1 & 4 & 1 & 0\\
        0 & 1 & 4 & 1\\
        0 & 0 & 1 & 2\\
    \end{smallmatrix}\right), \quad Q_1 = T_0E_1, \quad Q_2 = \rho_0 ^{-1}E_2,
\end{equation}
\begin{equation}
    J = \left( \begin{smallmatrix}
        0 & D \\ -D^\top & 0
    \end{smallmatrix}\right), \quad R = \left( \begin{smallmatrix}
        0 & 0 \\ 0 & 0
    \end{smallmatrix}\right),\quad B = \left( \begin{smallmatrix}
        B_1  \\ B_2
    \end{smallmatrix}\right),
\end{equation}
and 
\begin{equation}
    D = \tfrac{1}{2} \left(\begin{smallmatrix}
        -1 &  1 &  0 & 0\\
        -1 &  0 &  1 & 0\\
        0  & -1 &  0 & 1\\
        0  &  0 & -1 & 1\\
    \end{smallmatrix}\right), \quad B_1 =  \left(\begin{smallmatrix}
        0 &  0 \\
        0  &  0 \\
        0  &  0 \\
        0  &  0 \\
    \end{smallmatrix}\right), \quad B_2 =  \left(\begin{smallmatrix}
        -1 &  0 \\
        0  &  0 \\
        0  &  0 \\
        0  &  1 \\
    \end{smallmatrix}\right),
\end{equation}
for the case of Neumann boundary conditions \eqref{Neumann}, and 
\begin{equation}
    D = \tfrac{1}{2} \left(\begin{smallmatrix}
        1 &  1 &  0 & 0\\
        -1 &  0 &  1 & 0\\
        0  & -1 &  0 & 1\\
        0  &  0 & -1 & 1\\
    \end{smallmatrix}\right), \quad B_1 =  \left(\begin{smallmatrix}
        -1 &  0 \\
        0  &  0 \\
        0  &  0 \\
        0  &  0 \\
    \end{smallmatrix}\right), \quad B_2 =  \left(\begin{smallmatrix}
        0 &  0 \\
        0  &  0 \\
        0  &  0 \\
        0  &  1 \\
    \end{smallmatrix}\right)
\end{equation}
for the case of mixed boundary conditions \eqref{mixed}.
\end{example}

%\begin{remark}

\begin{remark}
\DM{Note that the application of PFEM to the Euler-Bernoulli beam has been presented in detail in \citep[\S 6.3, pp. 527--530]{Cardoso2020JournalPartitioned}; one of the features of 1D higher-order
operators is that many more variables are involved at the boundary; typically for the Euler-Bernoulli beam,
we end up with 4 boundary components instead of only 2 for the vibrating string: boundary values as well
as spatial derivatives of these boundary values have to be taken into account. Moreover on the numerical side, for the choice of the Finite Element families, since a higher regularity in space is required, we resort to Hermite polynomials instead of Lagrange polynomials. 
}
\end{remark}
\section{Model Order Reduction in the Loewner Framework}
\label{sec:Loewner}

In the previous section, we have presented the first approximation of the PDE \eqref{Eq:PDE}-\eqref{InputOutput} using the FEM, leading to the discretized model \eqref{StateSpace2}. For the sake of completeness of this paper and since the size of the discretized model can contain a large amount of states, in this section we recall the Loewner framework technique for the MOR. This approach has shown to be well suited both for large order systems and infinite-dimensional ones (see e.g. \citep{Poussot2023ConferenceData} and references therein). In this paper, we exploit the properties of this approach to generate a projector that allows to approximate the state of the discretized model in a structure preserving way.  

In the following, we rewrite the discretized port-Hamiltonian system as in the standard state space representation. To do this, we define $S = E^{-1}Q$, $A = (J-R)S$, and $C = B^\top S$ and the state-space representation in \eqref{StateSpace2} is equivalently written as
\begin{equation}\label{DiscretizedModel}
\Sigma_d :\begin{cases}
E\dot{x}_d (t) = Ax_d(t) + Bu(t), \quad x_d(0) = x_{d0},  \\
\quad \,\,  y(t) =  C x_d.
\end{cases}
\end{equation}
The objective of this section is to show the steps to build the matrices of a Reduced Order Model (ROM) of the following form:
\begin{equation} \label{Eq:ROM}
\Sigma_r : \begin{cases}
E_r \dot{x}_r(t) = A_r x_r(t) + B_r u(t), \quad 
x_r(0) = x_{r0},\\
\quad \,\,  y(t) = C_r x_r(t)
\end{cases}
\end{equation}
with $E_r,A_r \in \Real^{r\times r}$, $B_r \in \Real^{r \times n} $, $C_r \in \Real^{n \times r}$, and $x_{r0} \in \Real^{r}$, in which $r\ll N$, $n$ is the number of inputs and outputs of the initial PDE \eqref{Eq:PDE}-\eqref{InputOutput}, and $N = n_1N_1 + n_2N_2$ is the amount of states of the discretized model \eqref{DiscretizedModel}. Moreover, we provide a lifting projector $T$ that allows to approximate the state of the discretized (large-scale) model \eqref{DiscretizedModel} as follows:
\begin{equation}\label{Eq:T}
    x_d(t)\approx Tx_r(t),
\end{equation} recovering the physical meaning of the state variables.

\begin{remark}
    Note that since the discretized model \eqref{DiscretizedModel} has a port-Hamiltonian structure with quadratic storage function, it can be written in different ways as discussed in \citep{Willems1972Dissipative}, \citep{Beattie2022Port}. different structure-preserving techniques can be applied. Namely, interpolation techniques using {\it spectral zeros} \citep{Sorensen2005JournalPassivity}, \citep{Benner2020JournalIdentification}, \citep{Poussot2023ConferenceData}, or $\mathcal{H}_2$-inspired algorithms \citep{Gugercin2012JournalStructure}, for instance. See \citep[Section~8]{Mehrmann2023JournalControl} for a more complete discussion on structure-preserving MOR for Linear Time Invariant port-Hamiltonian Systems..
\end{remark}

\subsection{Standard solution in the Loewner approach}\label{SubSec:StandardLoewner}
We follow \citep[Section~3]{Antoulas2016ChapterATutorial} (see also \citep[Problem~5.1]{Mayo2007JournalFramework}) for the construction of the ROM. Let us define the transfer function associated to the discretized model \eqref{DiscretizedModel} as:
\begin{equation}\label{Eq:Transfer}
G(s) = C(sE-A)^{-1}B.
\end{equation}
This transfer function is used to generate the data used in the interpolation-based ROM construction. We define the right data as
\begin{equation}\label{rightData}
   (\lambda_j,r^\lambda_j,w_j), \quad w_j = G(\lambda_j)r_j ^\lambda, \quad j = \lbrace 1,\cdots , k \rbrace,
\end{equation}
with $\lambda_j \in \mathbb{C}$, $r_j^\lambda \in \mathbb{C}^{n \times 1}$, and $w_j \in \mathbb{C}^{n \times 1}$. Similarly, we define the left data as 
\begin{equation}\label{leftData}
    (\mu_i,l^\mu_i,v_i),\quad v_i = l_i^\mu G(\mu_i), \quad i = \lbrace 1,\cdots , q \rbrace,
\end{equation}
with $\mu_i \in \mathbb{C}$, $l_i^\mu \in \mathbb{C}^{1 \times n}$, and $v_j \in \mathbb{C}^{1 \times n}$. The objective is to construct the matrices $(E_r,A_r,B_r,C_r)$ by only using the right and left data \eqref{rightData}-\eqref{leftData}. Since the discretized model that we aim to reduce is passive, it has the same number of inputs than outputs $n$, implying the same sizes for the vectors $r_j^\lambda$ and $w_j$, and similarly with $l_i^\mu$ and $v_i$. For simplicity of the explanation, in the following we consider the same amount of left and right data and we denote it as $N_r = q = k$. 

Using the right and left data \eqref{rightData}-\eqref{leftData}, we build the {\it Loewner} ($\mathbb{L}$) and {\it shifted Loewner} ($\mathbb{\sigma L}$) matrices as follows:
\begin{equation}\label{Eq:Loewner}
\mathbb{L}_{ij} = \frac{v_i r^\lambda_j - {l^\mu_i}  w_j}{\mu_i - \lambda_j}, \quad\mathbb{\sigma L}_{ij} = \frac{\mu_i v_i r^\lambda_j - \lambda_j {l^\mu_i} w_j}{\mu_i - \lambda_j} .
\end{equation}
The main assumption concerning the right and left data \eqref{rightData}-\eqref{leftData}
is the following:
\begin{assumption}\label{Assumtion_k}
    For all $s \in \lbrace \lambda_j \rbrace  \cup \lbrace \mu_i \rbrace $; the following rank condition is satisfied:
    \begin{equation}
        rank(s\mathbb{L}-\sigma \mathbb{L}) = rank\left( \begin{bmatrix}
            \mathbb{L} & \sigma\mathbb{L}
        \end{bmatrix}\right)  = rank\left( \begin{bmatrix}
            \mathbb{L} \\ \sigma\mathbb{L}
        \end{bmatrix}\right)=k.
    \end{equation}
\end{assumption}
We arrange the right and left data \eqref{rightData}-\eqref{leftData} in the following matrix form:
\begin{equation}\label{LaMu}
\begin{split}
        &\Lambda = \begin{pmatrix}
        \lambda_1 & & 0 \\
         & \ddots & \\
         0 & & \lambda_{N_r}
    \end{pmatrix},  W = \begin{pmatrix}
w_1 & \cdots & w_{N_r}\end{pmatrix}, r = \begin{pmatrix}
r_1 & \cdots & r_{N_r}\end{pmatrix} , \\ 
    &M = \begin{pmatrix}
        \mu_1 & & 0 \\
         & \ddots & \\
         0 & & \mu_{N_r}
    \end{pmatrix},  V = \begin{pmatrix}
v_1 \\ \vdots \\ v_{N_r}\end{pmatrix}, \quad \quad \quad \quad  l = \begin{pmatrix}
l_1 \\ \vdots \\ l_{N_r}\end{pmatrix}.
\end{split}
\end{equation}
Then, a minimal realization using the right and left data \eqref{rightData}-\eqref{leftData} is given by:
\begin{equation} \label{Eq:1DWaveLoewner}
\Sigma_l : \begin{cases}
E_l \dot{x}_l(t) = A_l x_l(t) + B_l u(t), \quad 
x_l(0) = x_{l0},\\
\quad \,\,  y(t) = C_l x_l(t)
\end{cases}
\end{equation}
with $E_l = -\mathbb{L}$, $ A_l = -\mathbb{\sigma L}$, $B_l = V$, $C_l = W $. The transfer function:
\begin{equation}\label{Eq:TransferComplex}
G_l(s) = C_l(sE_l-A_l)^{-1}B_l.
\end{equation}
interpolates the data \eqref{rightData}-\eqref{leftData} throughout the tangential directions, {\it i.e.,} 
$G(\lambda_j)r_j^\lambda = G_l(\lambda_i)r_j^\lambda$ and $l_i^\mu G(\mu_j) = l_i \mu G_l(\mu_i)$ for all $j \in \lbrace 1,\cdots , N_r \rbrace$ and $i \in \lbrace 1,\cdots , N_r \rbrace$, respectively. Moreover, one can verify that the {\it tangential generalized controllability} ($C_b$) and { \it tangential generalized observability} ($O_b$) matrices, defined as 
\begin{equation}\label{ObCb}
\begin{split}
	&{O}_b = \begin{pmatrix}
	C(\mu_1E - A)^{-1} \\
	\vdots \\
	C(\mu_{N_r} E - A)^{-1}
	 \end{pmatrix}, \\
	& C_b = \begin{pmatrix}
	(\lambda_1E - A)^{-1}B &
	\cdots &
	(\lambda_{N_r} E - A)^{-1}B
	 \end{pmatrix}.
\end{split}
\end{equation}
satisfy $E_l = O_b E C_b$ and $A_l = O_b A C_b$. Then, $\Sigma_l$ can be seen as a transformation of $\Sigma_d$ though the {\it tangential generalized controllability} ($C_b$) and { \it tangential generalized observability} ($O_b$) matrices, {\it i.e.,}
if 
\begin{equation} \label{Eq:FirstProjection}
x_d = C_b x_l,
\end{equation}
one can compute the time derivative of \eqref{Eq:FirstProjection} multiplied by $E$ at the left side as follows:
\begin{align}
    E\dot{x}_d(t) &= EC_b \dot{x}_l(t), \\
    Ax_d(t)+B u(t) &= EC_b \dot{x}_l(t), \\
    AC_bx_l(t)+B u(t) &= EC_b \dot{x}_l(t), \\
    O_bAC_bx_l(t)+ O_bB u(t) &= O_bEC_b \dot{x}_l(t), \\
    A_lx_l(t)+ B_l u(t) &= E_l \dot{x}_l(t), 
\end{align}
where $O_bB = B_l$ can be verified by multiplying $O_b$ from \eqref{ObCb} by $B$ on the right and replacing $B_l$ by the definition in \eqref{leftData}-\eqref{LaMu}-\eqref{Eq:1DWaveLoewner}.

\subsection{Realization with real values:}\label{SubSec:RealValues}

Since we aim at finding a realization with matrices containing real entries, the interpolation points have to contain pairs of conjugated complex values, {\it i.e.,} the sets of right and left data \eqref{rightData}-\eqref{leftData} have to be selected in such a way that:
\begin{equation} \label{Eq:Conjugated}
\begin{split}
\lambda &=  \lbrace \lambda_1,\bar{\lambda}_2, \cdots,\lambda_{M_r}, \bar{\lambda}_{M_r} \rbrace, \quad 
\mu =  \lbrace \mu_1,\bar{\mu}_2, \cdots,\mu_{M_r}, \bar{\mu}_{M_r} \rbrace,
\end{split}
\end{equation}
with the 'bar' notation indicating the conjugated value and $2M_r = N_r$. We notice that by incorporating the conjugated values, the amount of data is even. One could also add strictly real interpolation points, but for simplicity of the explanation, we use only complex values. The realization \eqref{Eq:1DWaveLoewner} is composed then by complex-valued matrices, {\it i.e.,} $E_l$, $A_l \in \mathbb{C}^{N_r \times N_r}$, $B_l \in \mathbb{C}^{N_r \times n}$, and $C_l \in \mathbb{C}^{n \times N_r}$. Then, using the data as in \eqref{Eq:Conjugated}, \eqref{Eq:1DWaveLoewner} can be transformed into a realization with real entries using the following projection:
\begin{equation} \label{Eq:SecondProjection}
x_l = \mathcal{J}\bar{x}_l
\end{equation}
 with $\mathcal{J} \in \mathbb{C}^{N_r \times N_r}$ such that
\begin{equation}
\mathcal{J} = \begin{pmatrix}
\mathcal{J}_i & \cdots & 0\\
\vdots & \ddots & \vdots \\
0 & \cdots & \mathcal{J}_i
\end{pmatrix}, \quad \mathcal{J}_i = \frac{1}{\sqrt{2}}\begin{pmatrix}
1 & -i \\ 1 & i
\end{pmatrix}.
\end{equation}
Then, the realization with real matrices is written as:
\begin{equation} \label{Eq:1DWaveLoewner_real}
\bar{\Sigma}_l : \begin{cases}
\bar{E}_l \dot{\bar{x}}_l(t) = \bar{A}_l \bar{x}_l(t) + \bar{B}_l u(t), \quad 
\bar{x}_l(0) = \bar{x}_{l0},\\
\quad \,\,  y(t) = \bar{C}_l \bar{x}_l(t)
\end{cases}
\end{equation}
with $\bar{E}_l = \mathcal{J}^* E_l \mathcal{J}$, $\bar{A}_l = \mathcal{J}^* A_l \mathcal{J}$, $\bar{B}_l = \mathcal{J}^* B_l$, $\bar{C}_l =  C_l \mathcal{J}$. One can guarantee that the transfer matrix
\begin{equation}\label{Eq:TransferReal}
\bar{G}_l(s) = \bar{C}_l(s\bar{E}_l-\bar{A}_l)^{-1}\bar{B}_l, 
\end{equation}
interpolates the data \eqref{rightData}-\eqref{leftData} throughout the tangential directions, {\it i.e.,} 
$G(\lambda_j)r_j^\lambda = \bar{G}_l(\lambda_i)r_j^\lambda$ and $l_i^\mu G(\mu_j) = l_i \mu \bar{G}_l(\mu_i)$ for all $j \in \lbrace 1,\cdots , N_r \rbrace$ and $i \in \lbrace 1,\cdots , N_r \rbrace$, respectively. 

\subsection{Order reduction via Singular Value Decomposition}\label{SubSec:SVD}
One of the main properties of the {\it Loewner} framework is the 'rank revealing' characteristic \citep{Gosea2015ConferenceModel}. This means that if we have enough data $N_r$, the rank of the {\it Loewner} matrix will reveal the McMillan degree $k$ required for the ROM among the selected data (see \citep[Lemma~2.1]{Mayo2007JournalFramework}). Since the selected data \eqref{rightData}-\eqref{leftData} may contain redundant values, one can build the following short Singular Value Decompositions (SVDs):
\begin{equation}
    \begin{bmatrix}
        \bar{E}_l & \bar{A}_l
    \end{bmatrix} = Y\Sigma_l \tilde{X}^*, \quad \quad \begin{bmatrix}
        \bar{E}_l \\ \bar{A}_l
    \end{bmatrix} = \tilde{Y}\Sigma_r X^*, 
\end{equation}
with $Y$, $X \in \mathbb{C} ^{N_r \times k}$, $\Sigma_l$, $\Sigma_r \in \mathbb{C}^{k \times k}$, and $k$ satisfying Assumption \ref{Assumtion_k}. Define $E_r =Y^*\bar{E}_l X$, $A_r =Y^*\bar{A}_l X$, $B_r =Y^*\bar{B}_l$, and $C_r =\bar{C}_l X$.  Then, the realization \eqref{Eq:ROM} with transfer function
\begin{equation}\label{Eq:TransferROM}
{G}_r(s) = {C}_r(s{E}_r-{A}_r)^{-1}{B}_r, 
\end{equation}is a descriptor realization of an approximate interpolant of the data with McMillan degree $k = rank(\mathbb{L})$ (see \citep[Theorem~2.19]{Antoulas2016ChapterATutorial}). Similarly as in \eqref{Eq:FirstProjection}, one can approach the state $\bar{x}_l$ of \eqref{Eq:1DWaveLoewner_real} through the following projection
\begin{equation}\label{Eq:LastProj}
    \bar{x}_l(t) \approx X x_r(t).
\end{equation}
Thus, by combining the projection \eqref{Eq:FirstProjection}, \eqref{Eq:SecondProjection}, and \eqref{Eq:LastProj}, one can build the projector $T$ of \eqref{Eq:T} as follows:
\begin{equation}\label{Eq:TDef}
    T := C_b\mathcal{J}X.
\end{equation}

\subsection{Discussion about the structure-preserving ROM}

Until here, the Loewner approach offers an efficient algorithm to construct a ROM for the discretized model that has been obtained using the FEM. We emphasize that the main numerical computational steps are:
\begin{itemize}
    \item the inversion of the matrices $(sE-A)$ for all $s \in \lbrace \lambda_j \rbrace  \cup \lbrace \mu_i \rbrace$ in order to build the right and left data \eqref{rightData}-\eqref{leftData} and the {\it tangential generalized controllability} ($C_b$) and { \it tangential generalized observability} ($O_b$) matrices in \eqref{ObCb},
    \item the construction of the {\it Loewner} and {\it shifted Loewner} matrices in \eqref{Eq:Loewner}, which after building the right and left data \eqref{rightData}-\eqref{leftData}, the numerical efforts are negligible.
\end{itemize}
However, the obtained ROM \eqref{Eq:ROM} is not always guaranteed to be in a port-Hamiltonian form. In particular, using the classic approach \citep{Antoulas2016ChapterATutorial}, the interpolation points are chosen along the imaginary axis ($s_i= \sqrt{-1} \omega _i$) for some desired frequencies $\omega_i$ and tangential directions $r_i$ as unitary vectors. One could iterate the construction of the realization \eqref{Eq:ROM} by modifying the frequencies $\omega_i$ until a ROM is obtained preserving the port-Hamiltonian structure.

As it is discussed in \citep{Antoulas2008JournalConstruction}, by properly choosing the interpolation points $s_i$ and tangential directions $r_i$, one can force the passivity of the reduced-order model. In \citep{Benner2020JournalIdentification} and \citep{Poussot2023ConferenceData}, the authors have proposed algorithms that lead to port-Hamiltonian realizations by selecting the frequency-domain data $s_i$ as the {\it spectral zeros} and the tangential direction $r_i$ as the {\it zero directions}. In \citep{Poussot2023ConferenceData}, for the cases in which $D=0$, the authors have proposed to shift the data with a feedthrough term $D_r$ and shifting back the obtained realization using the same $D_r$.

Following \citep[Algorithm~1]{Benner2020JournalIdentification} and \citep[Algorithm~1]{Poussot2023ConferenceData}, we define the {\it spectral zeros} $s_i$ and {\it zero directions} $r_i$ associated to the realization \eqref{Eq:ROM} shifted by $D=D_r$ with $D_r + D_r^\top > 0$ as the solutions of the following generalized eigenvalue problem:
\begin{equation}\label{Eq:SpectralZeros}
\begin{bmatrix}
0 & A_r & B_r \\
A_r^\top & 0 & C_r^\top \\
B_r^\top & C_r & D_r+D_r^\top
\end{bmatrix}\begin{bmatrix}
p_i \\ q_i \\ r_i
\end{bmatrix} = s_i\begin{bmatrix}
0 & E_r & 0 \\
-E_r^\top & 0 & 0 \\
0 & 0 & 0
\end{bmatrix}\begin{bmatrix}
p_i \\ q_i \\ r_i
\end{bmatrix}.
\end{equation}
From all solutions of \eqref{Eq:SpectralZeros}, we select all the values $(s_i,r_i)$ such that the real part of $s_i$ is such that $0<Re(s_i)<\infty$. Then, we re-define the right and left data as in \eqref{rightData}-\eqref{leftData}, using the following values:
\begin{equation}\label{NewData}
\begin{split}
     &\lambda_i = s_i, \quad  \;\; r_i^\lambda = r_i, \quad w_i = \left( G_r(\lambda_i) + D_r\right)r_i ^\lambda, \\
     &\mu_i = -\bar{s}_i, \quad l_i^\mu = r_i^*,\quad v_i = l_i^\mu \left( G_r(\lambda_i) + D_r\right) \\
\end{split}
\end{equation}
with $G_r$ defined in \eqref{Eq:TransferROM} and $D_r$ free to choose. Then, using the new data \eqref{NewData}, we repeat again the computations presented in Subsection \ref{SubSec:StandardLoewner} and Subsection \ref{SubSec:RealValues} (the SVD of Subsection \ref{SubSec:SVD} is not required at this stage). Finally, we remark that since the Loewner matrix in \eqref{Eq:Loewner} is positive definite, an explicit port-Hamiltonian realization can be written from the ROM. See \citep[Alorithm~1,~Step~6]{Benner2020JournalIdentification} for further details.
\section{Summary of the procedure}\label{Sec:Summary}
%\textbf{Input:} The matrices \eqref{PDEMatrices} of the PDE \eqref{Eq:PDE} and the input and output matrices \eqref{InputOutput}. 

In the following, we summarize the complete procedure for the discretization and reduction of the initial PDE \eqref{Eq:PDE}-\eqref{InputOutput}. 
\begin{enumerate}
\item Discretize the PDE \eqref{Eq:PDE}-\eqref{InputOutput} using the PFEM proposed in Section \ref{sec:PFEM}. This leads to the port-Hamiltonian realization \eqref{StateSpace2} that can be equivalently written as in \eqref{DiscretizedModel}.
\item Follow Subsection \ref{SubSec:StandardLoewner}, \ref{SubSec:RealValues} and \ref{SubSec:SVD} to find a preliminary ROM \eqref{Eq:ROM} using the standard Loewner approach, in which the interpolation points are on the imaginary axis and satisfy Assumption \ref{Assumtion_k}.
\item If the obtained model \eqref{Eq:ROM} is passive, one can compute the projector \eqref{Eq:TDef} and \eqref{Eq:ROM} can be used to approximate the discretized model \eqref{StateSpace2} using the projection \eqref{Eq:T}.
\item If the preliminary ROM \eqref{Eq:ROM} is not passive, then compute the {\it spectral zeros} $s_i$ and {\it zero directions} $r_i$ solving the generalized eigenvalue problem in \eqref{Eq:SpectralZeros}, using a shift term $D_r$ such that $D_r+D_r^\top>0$.
\item Redefine the right and left data as in \eqref{NewData} and repeat Subsection \ref{SubSec:StandardLoewner} and \ref{SubSec:RealValues}. The ROM is as the one defined in \eqref{Eq:1DWaveLoewner_real} and the projection is given by $x_d(t) \approx T \bar{x}_l(t)$ with $T =C_b \mathcal{J}$. 
\end{enumerate}

\section{Application cases}
\label{Examples}

In this section, the one-dimensional wave equation is firstly used to exemplify the discretization and reduction methodology proposed in this paper. Secondly, a slightly more complex system, namely the Timoshenko beam is also used to show the versatility of the proposed approach when changing the PDE and the boundary inputs and outputs.

\begin{example}\label{Example3}
(One-dimensional wave equation) In the following, we consider the case of mixed boundary conditions studied in Example \ref{Example2}, with unitary parameters ($T_0=\rho_0 = 1$, $a = 0$, and $b = 1$), and a higher order of discretization using $N=500$ basis functions, which gives a total of  $nN = 1000$ states. First, we choose interpolation points that lie on the imaginary axis in a frequency range $\omega \in [0.9,8.5]$ $rad/s$. and the tangential directions are interspersed between $\left(\begin{smallmatrix}
    1 \\ 0 
\end{smallmatrix}\right)$ and $\left(\begin{smallmatrix}
    0 \\ 1 
\end{smallmatrix}\right)$. In Figure \ref{Fig:BodeG22}, we show the magnitude of the bode diagram of the transfer function relating the second input and second output (force and velocity at the right, respectively). In solid blue, we show the discretized model frequency response, in dashed orange, the first ROM obtained using the standard version of the Loewner approach, and in dashed green the second ROM using the spectral zeros and zero directions as interpolation points. The right and left data are shown in dots violet and yellow, respectively. As  expected, both ROMs are closer to the discretized model near the frequencies range that has been selected. We notice that on the one hand, the preliminary ROM matches the $\lbrace \lambda_i \rbrace \cup \lbrace \mu_i \rbrace$ due to the fact that the interpolations points lies on the imaginary axis and the tangential directions are unitary vectors (as mentioned above). On the other hand, the final ROM, which uses spectral zeros and zero directions as interpolation points and tangential directions, respectively, matches less well than the preliminary ROM. However, this is the price to pay in order to obtain a passive ROM. This important constrain for the final ROM allows us to preserve stability. 
\begin{figure}[!h]
\begin{center}
\includegraphics[width=0.48\textwidth]{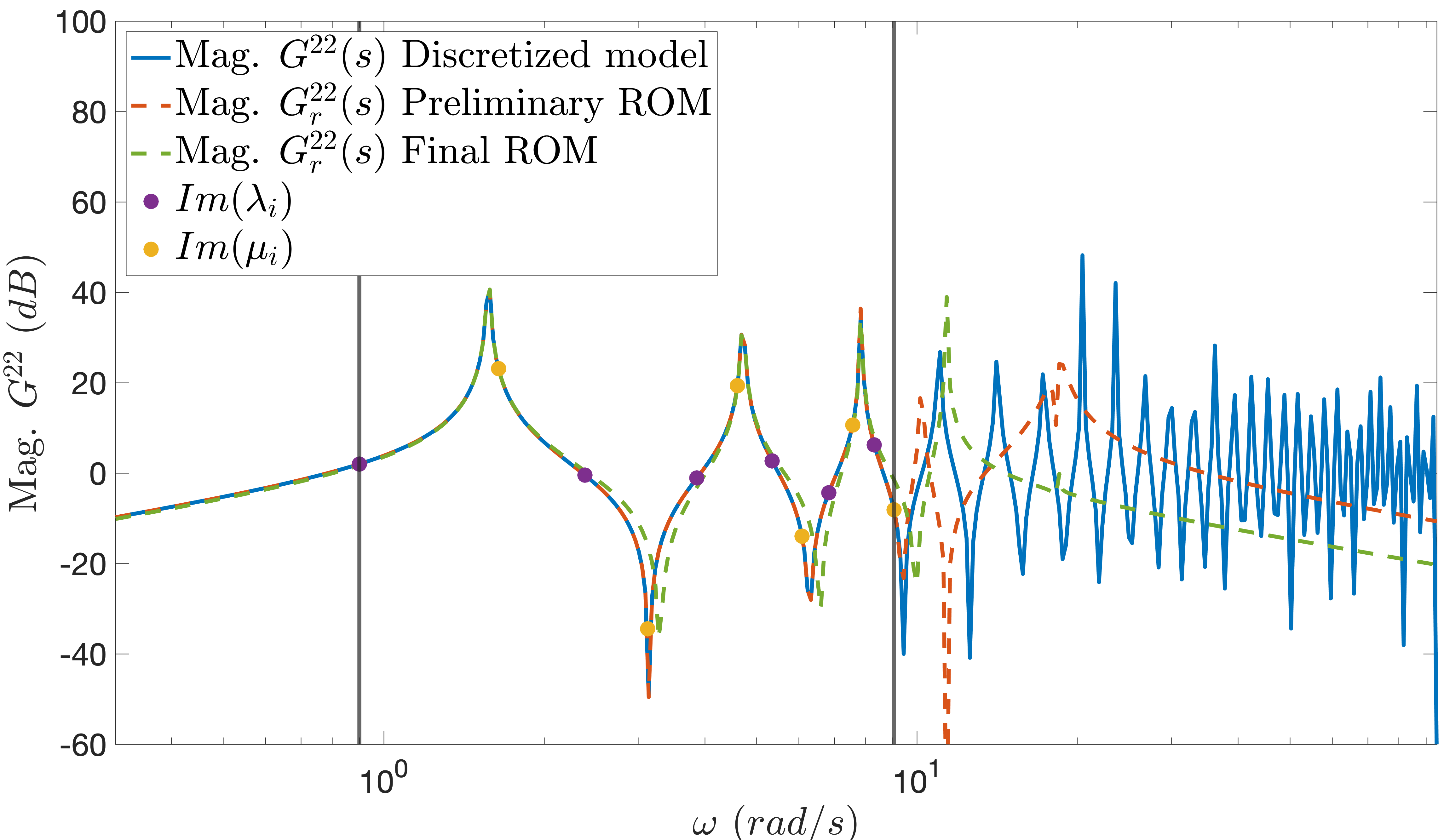}  
\caption{Magnitude Bode diagrams of the discretized model and the reduced one.} 
\label{Fig:BodeG22}    
\end{center}
\end{figure}

In Figure \ref{Fig:Hamiltonians1}, we show the energies obtained from a numerical simulation of the discretized model, the preliminary ROM, and the final ROM. For all simulations, the string is considered to be attached at the left side. We consider null initial condition and we exite the system with a force at the right side with the form $u(t) = sin(1.6t)$ for $t\in[0,5]$ and $u(t) =0$ for $t>5$. We can see that the preliminary version of the ROM diverges when $u = 0$ and that the final ROM is conservative, {\it i.e.,} $\dot{H} = 0$.
\begin{figure}[!h]
\begin{center}
\includegraphics[width=0.48\textwidth]{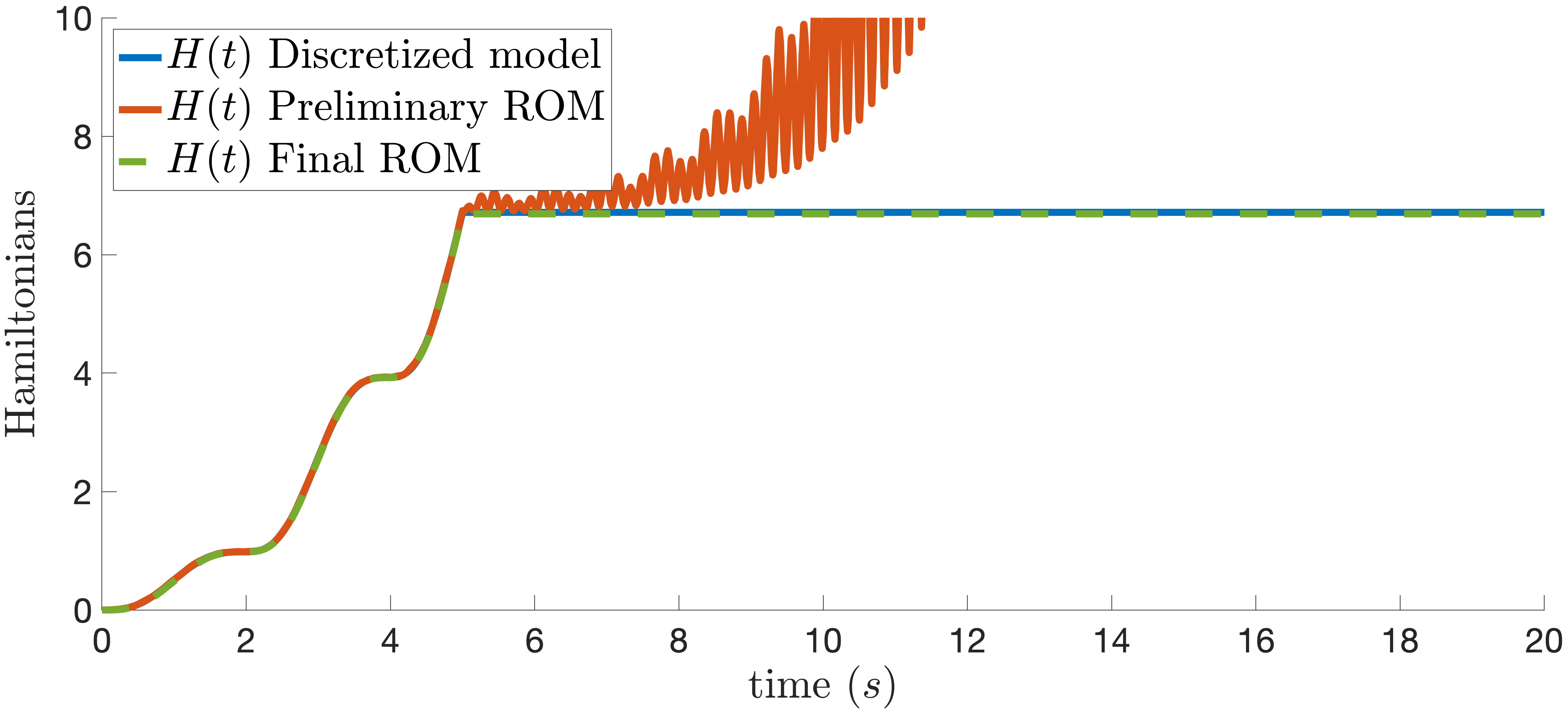}  
\caption{Energy of the discretized model, preliminary ROM and final ROM.} 
\label{Fig:Hamiltonians1}    
\end{center}
\end{figure}
\end{example}

\begin{example}\label{sec:TBExample}(Timoshenko Beam)
We consider the following Timoshenko beam model:
\begin{equation}\label{TimoshenkoBeam}
\begin{split}
\rho(\zeta) \frac{\partial ^2 w}{\partial t ^2}(\zeta,t) =& \frac{\partial}{\partial \zeta}\left( K(\zeta)\left(\frac{\partial w}{\partial \zeta} (\zeta,t) - \theta(\zeta,t) \right) \right) - g_1(\zeta)\frac{\partial w}{\partial t} (\zeta,t), \\
I_{\rho}(\zeta) \frac{\partial ^2 \theta}{\partial t ^2}(\zeta,t) =& \frac{\partial}{\partial \zeta}\left( EI(\zeta)\frac{\partial \theta}{\partial \zeta} (\zeta,t)  \right)- g_2(\zeta)\frac{\partial \theta}{\partial t} (\zeta,t) \\ 
&\quad \quad + K(\zeta) \left( \frac{\partial w}{\partial \zeta}(\zeta,t) - \theta(\zeta,t) \right).
\end{split}
\end{equation}
where $\zeta \in [a,b]$ is the spatial domain, $t \geq 0$ is the time, $w(\zeta,t)$ is the transverse displacement of the beam, and $\theta(\zeta,t)$ is the rotation angle of a filament of the beam. The space-dependent energy parameters are the shear modulus $K(\zeta)$, the Young's modulus of elasticity multiplied by the moment of inertia of the cross section $EI(\zeta)$, the mass density $\rho(\zeta)$, and the rotatory moment of inertia of a cross section $I_\rho(\zeta)$. $K$, $EI$, $\rho$ and $I_\rho$ are positive definite functions for all $\zeta \in [a,b]$. The in-domain  dissipation parameters are $g_1(\zeta)$ and $g_2(\zeta)$. Both are semi-positive definite functions for all $\zeta \in [a,b]$. 
The Hamiltonian of \eqref{TimoshenkoBeam} is 
\begin{equation}\label{HTimoshenkoBeam}
\begin{split}
H=\frac{1}{2}\int_a^b K & \left(\frac{\partial w}{\partial \zeta}-\theta \right)^2 + EI\left(\frac{\partial \theta}{\partial \zeta} \right)^2 + {\rho}\left(\frac{\partial w}{\partial t} \right)^2+ {I_\rho}\left(\frac{\partial \theta}{\partial t} \right)^2 d\zeta,
\end{split}
\end{equation}
where the time and spatial dependency has been omitted for simplicity. We choose as state variables $x_1(\zeta,t) = \dfrac{\partial w}{\partial \zeta} (\zeta,t) - \theta(\zeta,t)$, $x_2(\zeta,t) = \dfrac{\partial \theta}{\partial \zeta} (\zeta,t)$, $x_3(\zeta,t) = \rho(\zeta)\dfrac{\partial w}{\partial t} (\zeta,t)$, and $x_4(\zeta,t) = I_\rho(\zeta)\dfrac{\partial \theta}{\partial t} (\zeta,t)$ corresponding to the shear displacement, the angular displacement, the momentum and the angular momentum, respectively. Then, equations \eqref{TimoshenkoBeam} can be written as:
\begin{equation}
\dfrac{\partial}{\partial t} \begin{pmatrix}
x_1 \\ x_2 \\ x_3 \\ x_4
\end{pmatrix} = \begin{pmatrix}
0 & 0 & \tfrac{\partial}{\partial \zeta} & -1 \\
0 & 0 & 0 & \tfrac{\partial}{\partial \zeta} \\
\tfrac{\partial}{\partial \zeta} & 0 & -g_1 &0  \\
1 &\tfrac{\partial}{\partial \zeta} & 0 & -g_2  \\
\end{pmatrix}\begin{pmatrix}
e_1 \\ e_2 \\ e_3 \\ e_4
\end{pmatrix}, \quad \begin{pmatrix}
    e_1  \\
    e_2  \\
    e_3 \\
    e_4
\end{pmatrix} = \begin{pmatrix}
    Kx_1 \\
    EIx_2 \\
    \frac{1}{\rho}x_3, \\
    \frac{1}{I_\rho}x_4, 
\end{pmatrix}.
\end{equation}
in which $e_1$, $e_2$, $e_3$, and $e_4$ correspond to the force, torque, transversal velocity, and angular velocity, respectively. We notice that the PDE takes the form of \eqref{Eq:PDE} with $P_{11} = P_{22} =G_{11}= \left(\begin{smallmatrix}
    0 & 0 \\ 0 & 0
\end{smallmatrix}\right)$, $P_{12} = P_{21} = \left(\begin{smallmatrix}
    1 & 0 \\ 0 & 1
\end{smallmatrix}\right)$, $G_{22} = \left(\begin{smallmatrix}
    g_1 & 0 \\ 0 & g_2
\end{smallmatrix}\right)$, $G_{12} = \left(\begin{smallmatrix}
    0 & 1 \\ 0 & 0
\end{smallmatrix}\right)$, $G_{21} = -G_{12}^\top$, $\mathcal{H}_1 = \left(\begin{smallmatrix}
K & 0 \\ 0 & EI
\end{smallmatrix}\right)$, and $\mathcal{H}_2 = \left(\begin{smallmatrix}
\rho^{-1} & 0 \\ 0 & I_\rho ^{-1}
\end{smallmatrix}\right)$. We consider the case in which the inputs are the velocities at the left side and forces at the right side with conjugated outputs. This corresponds to the following inputs and outputs:
\begin{equation}
u(t) = \begin{pmatrix}
e_3(a,t)\\
e_4(a,t)\\
e_1(b,t)\\
e_2(b,t)
\end{pmatrix},  \quad y(t) = \begin{pmatrix}
-e_1(a,t)\\
-e_2(a,t)\\
e_3(b,t)\\
e_4(b,t)
\end{pmatrix},
\end{equation}
\begin{equation}
\begin{bmatrix} V_{\mathcal{B}} \\ \hdashline[2pt/2pt]
V_{\mathcal{C}}
\end{bmatrix}
 = \begin{bmatrix}
0 & 0 & 0 & 0 & 0 & 0 & {1} & {0} \\
0 & 0 & 0 & 0  & 0 & 0 & {0} & {1} \\ 
{1} & {0}& 0 & 0 & 0 & 0  & 0 & 0 \\
{0} & {1} & 0 & 0 & 0 & 0 & 0 & 0 \\ \hdashline[2pt/2pt]
0 & 0 & 0 & 0 & -1 & 0 & 0 & 0 \\
0 & 0 & 0 & 0 & 0 & -1 & 0 & 0 \\
0 & 0 & 1 & 0 & 0 & 0 & 0 & 0 \\
0 & 0 & 0 & 1 & 0 & 0 & 0 & 0 \\
\end{bmatrix}.
\end{equation}

In the following, We consider the beam clamped at the left side and with a force and torque actuators at the right side, and collocated velocity and angular velocity sensors. For simplicity, we consider the following parameters $K = EI = \rho = I_\rho = 1$, $a = 0$ and $b = 1$, $g_1 = g_0 = 0$.

For the spatial discretization, we consider $\phi_1 = \phi_2 =: \phi$, and since $n_1=2$, and $n_2 = 2$ then, $\Phi_1 = \Phi_2 =:\Phi = \left( \begin{smallmatrix}
\phi & 0 \\ 
0 & \phi
\end{smallmatrix}\right)$. 
Then, one can compute the matrices $E$, $J$, $R$, $Q$, and $B$ of the discretized model \eqref{StateSpace2}, which leads to similar matrices as the ones obtained in Example \ref{Example2}. We emphasize that for different parametrizations of the input and output matrices $V_{\mathcal{B}}$ and $V_{\mathcal{C}}$, the construction of the matrix $J$ in \eqref{Eq:JE} can be easily adapted. In the following, we consider $N=500$ basis functions and since $n_1=n_2 = 2$, the discretized model contains $(n_1+n_2)N = 2000$ states.

For the MOR, we select frequencies $\omega \in [0.1, 20] \,(rad/s)$. The final ROM is obtained with $k = 32$ states. In Fig. \ref{fig:BodeBeam}, we show the frequency responses of the transfer function of the discretized model and the final ROM relating the third and fourth inputs and outputs. As expected, for frequencies greater than $\omega =20  \,(rad/s)$, the final ROM is no longer close to the discretized model.
\begin{figure}[!h]
     \centering
     \begin{subfigure}[b]{0.23\textwidth}
         \centering
         \includegraphics[width=\textwidth]{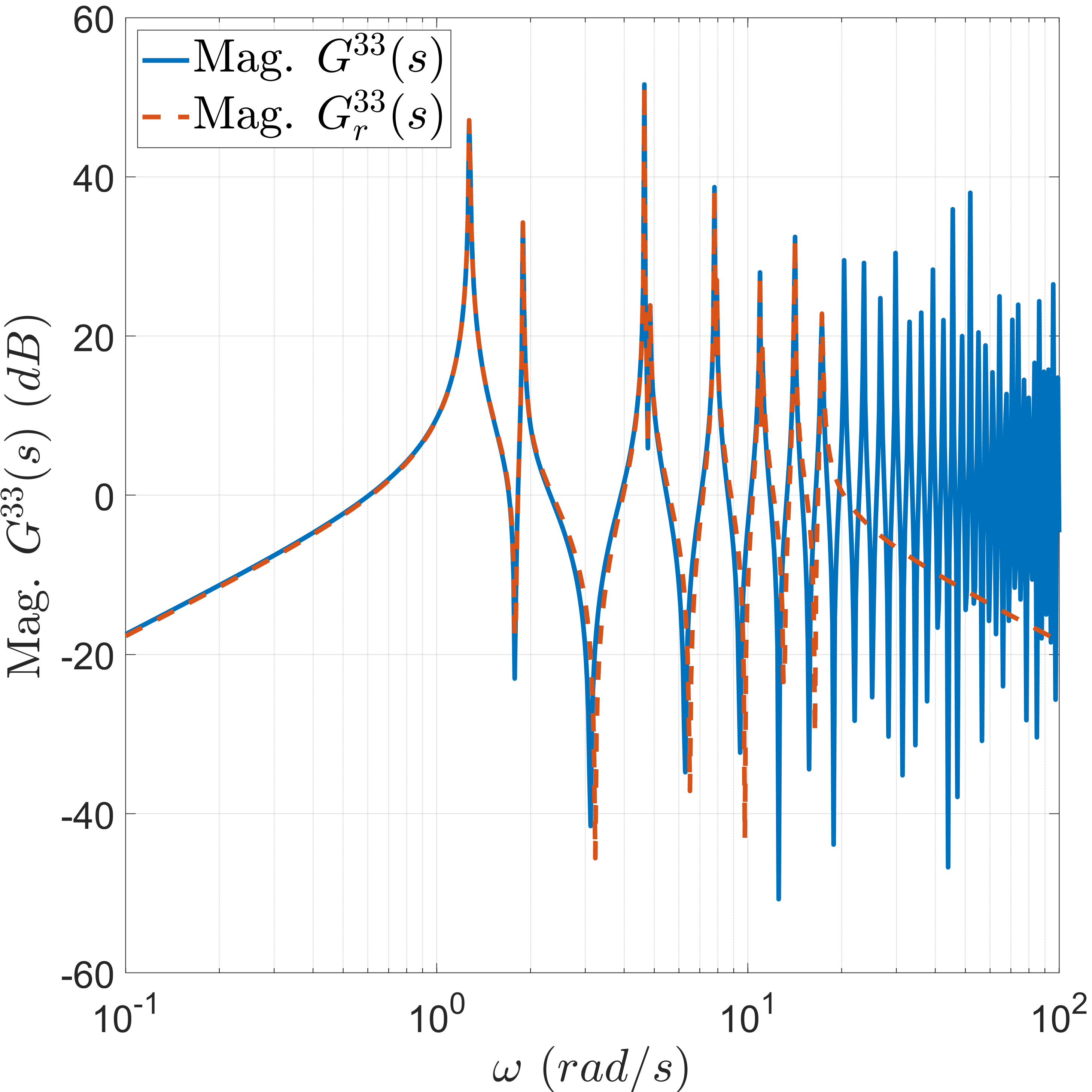}
%         \caption{}
%         \label{fig:y equals x}
     \end{subfigure}
     \hfill
     \begin{subfigure}[b]{0.23\textwidth}
         \centering
         \includegraphics[width=\textwidth]{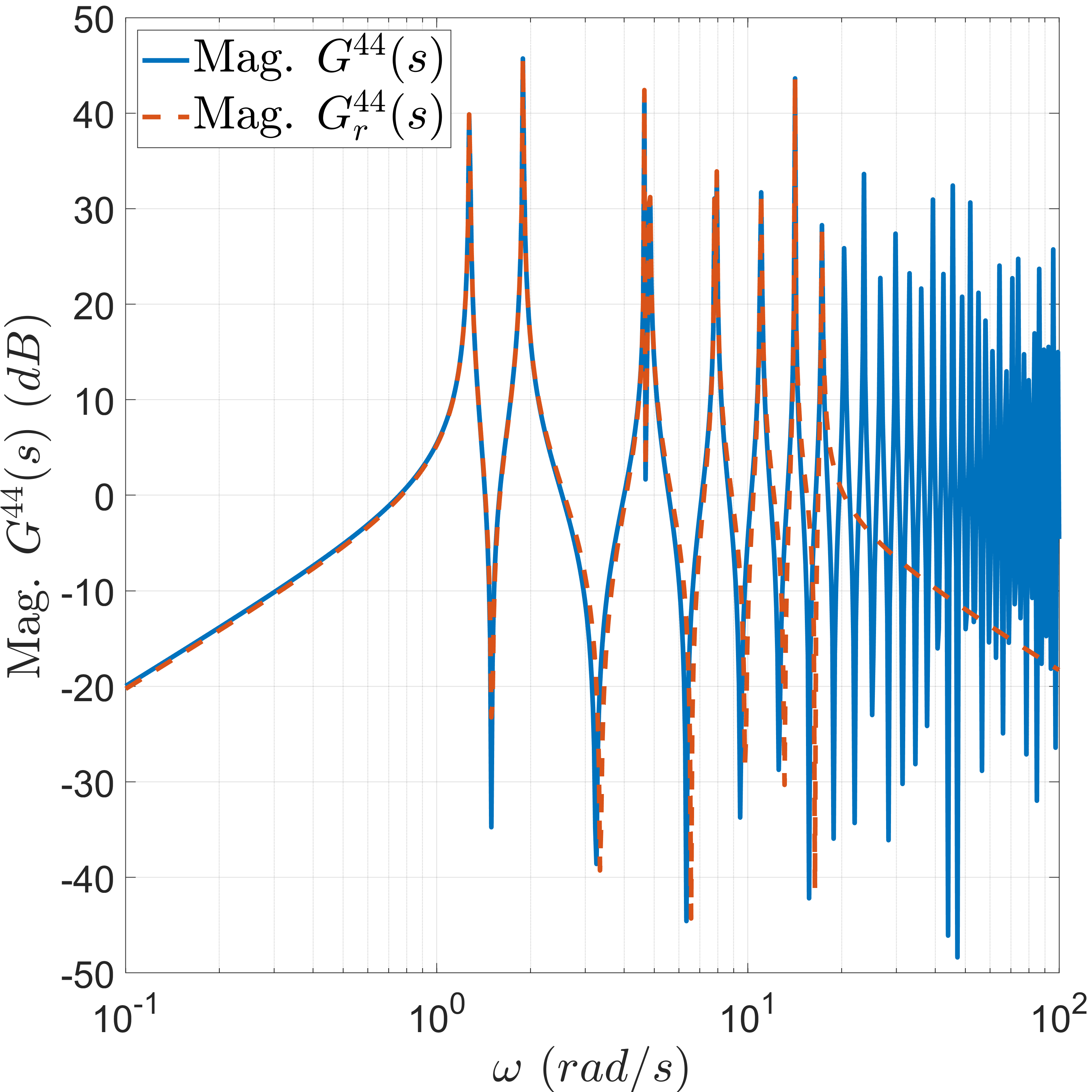}
%         \caption{}
%         \label{fig:three sin x}
     \end{subfigure}
        \caption{Bode magnitude diagrams of the discretized Timoshenko beam model and the ROM.}
        \label{fig:BodeBeam}
\end{figure}

Finally, we show a time simulation of the discretized model and the ROM using a step time $\delta t = 10^{-4}\,(s)$ and initial data $x_1(\zeta,0)=x_2(\zeta,0)=x_3(\zeta,0)=x_4(\zeta,0) = 0$. The simulation considers the beam attached at the left side with force actuators and velocity sensors at the right side. In this case, we simulate both models in closed-loop $(CL)$ using $u(t) =-Ky(t)+r$ with $K = diag(0,0,0.1,0.1)$ and $r = (0,0,1,-2)$. The output feedback gain $K$ represents a small damping at the right boundary through the force actuators and the signal $r$ is a change of equilibrium when applying a force of $1\,(N)$ and torque of $-2\,(Nm)$ at the right side of the beam. In Fig. \ref{Fig:HamiltonianBeam}, we show the Hamiltonian of both models. The travelling signals, natural from hyperbolic PDEs, can be appreciated in the Hamiltonian at every second. One can see that the ROM also capture this behaviour.  
%The closed-loop eigenvalues when using the output feedback are shown in Fig. \ref{fig:EigenvaluesBeam}. All the eigenvalues have negative real part and one can see that some of the low-frequencies values of the ROM approaches the ones of the PFEM.  
\begin{figure}[!h]
\begin{center}
\includegraphics[width=0.48\textwidth]{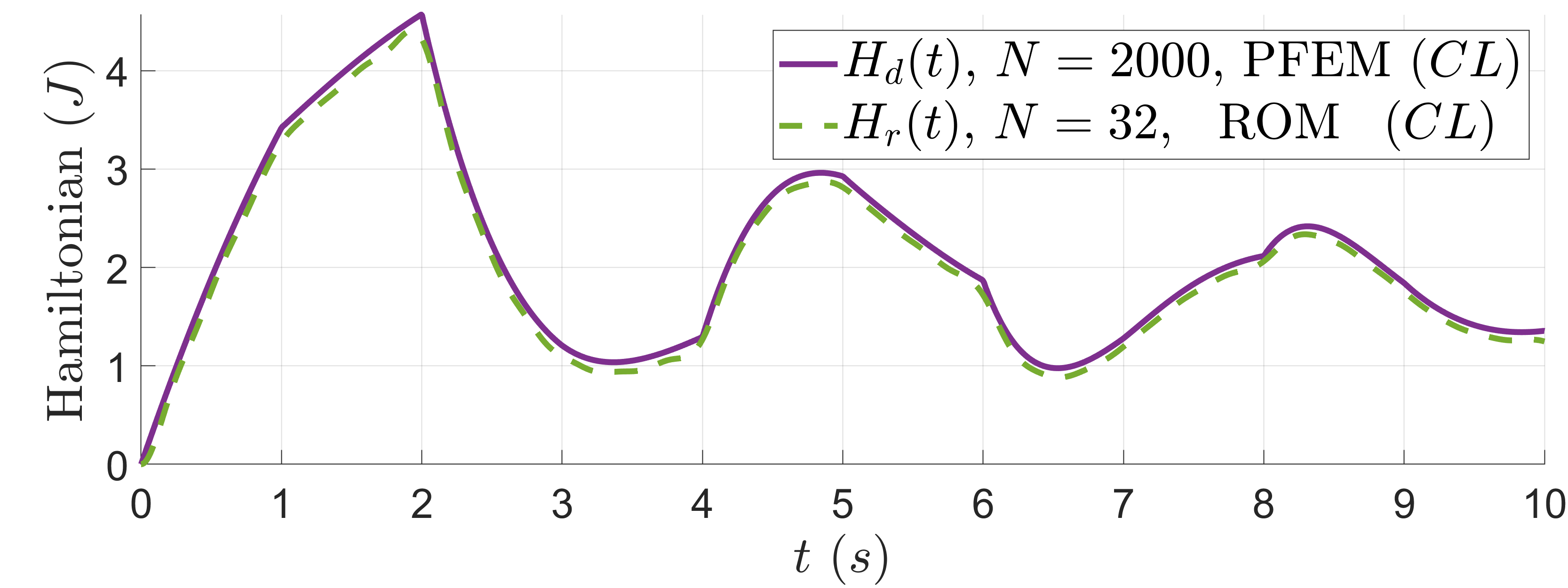}  
\caption{Hamiltonian of the PFEM model, and  the ROM in closed-loop $(CL)$ with $u = -K\,y+r$.} 
\label{Fig:HamiltonianBeam}                           
\end{center}
\end{figure}

Finally, the beam deformation is shown at different time screenshots in Fig. \ref{fig:SimulationBeam}. One can see that, in this case, since low-frequencies are predomiant in  the behaviour, the ROM is able to precisely simulate the discretized model. The complete simulations can be found in the github repository\footnote{\href{https://gitlab.com/ToledoZucco/drb}{https://gitlab.com/ToledoZucco/drb}}.
\begin{figure}[!h]
     \centering
     \begin{subfigure}[b]{0.48\textwidth}
         \centering
         \includegraphics[width=\textwidth]{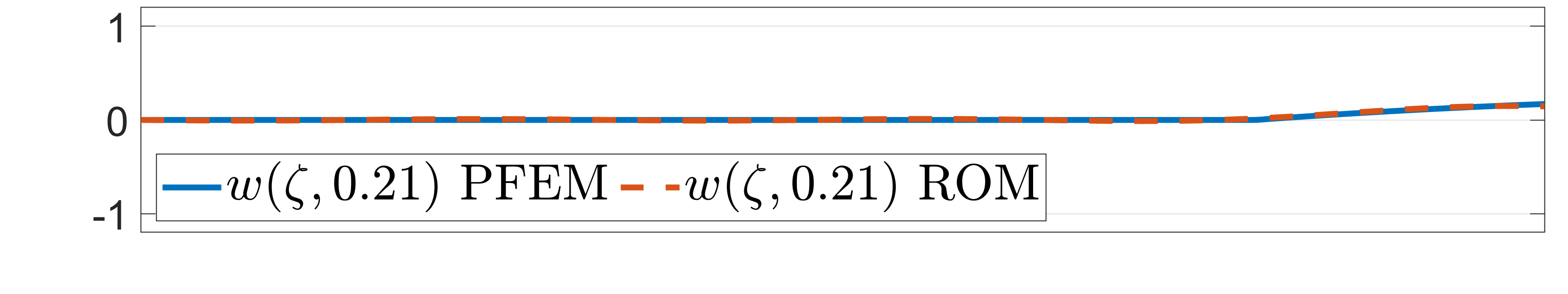}
%         \caption{}
%         \label{fig:SpectralZerosWave_a}
     \end{subfigure}
     \hfill
     \begin{subfigure}[b]{0.48\textwidth}
         \centering
         \includegraphics[width=\textwidth]{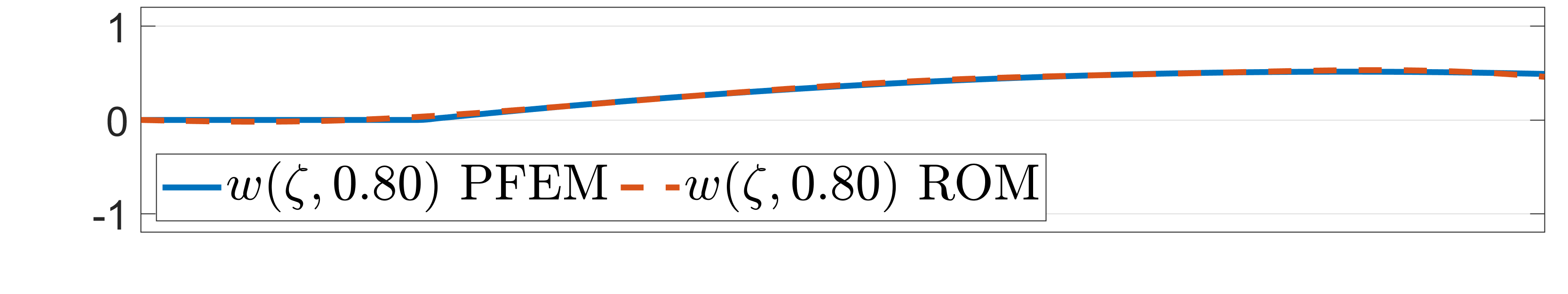}
%         \caption{}
%         \label{fig:SpectralZerosWave_b}
     \end{subfigure}
     \hfill
     \begin{subfigure}[b]{0.48\textwidth}
         \centering
         \includegraphics[width=\textwidth]{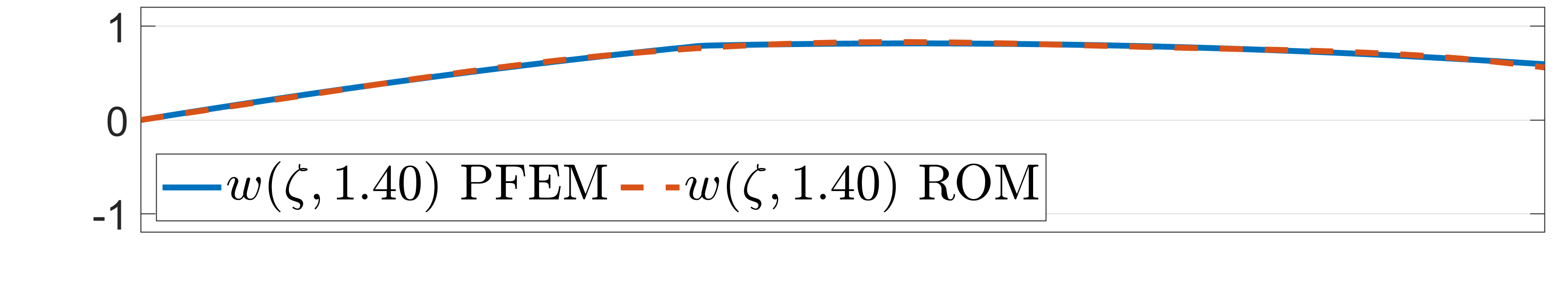}
%         \caption{}
%         \label{fig:SpectralZerosWave_b}
     \end{subfigure}
     \hfill
     \begin{subfigure}[b]{0.48\textwidth}
         \centering
         \includegraphics[width=\textwidth]{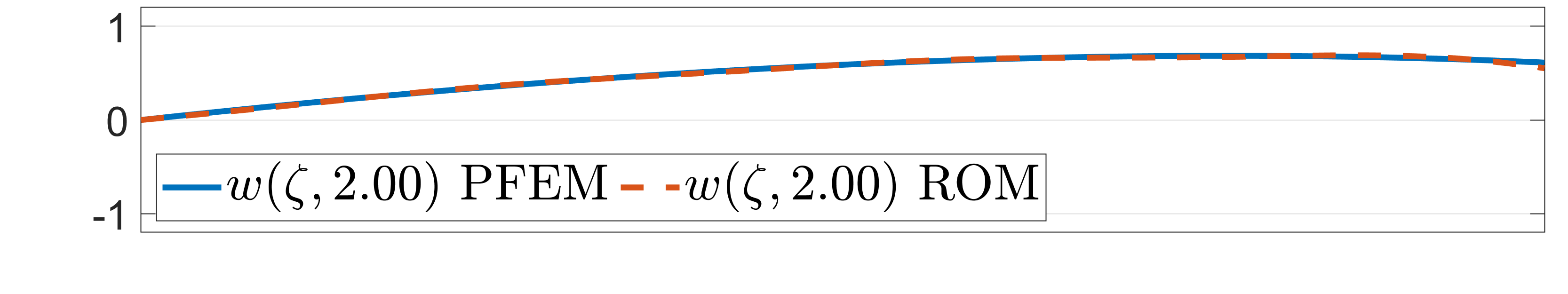}
%         \caption{}
%         \label{fig:SpectralZerosWave_b}
     \end{subfigure}
     \hfill
     \begin{subfigure}[b]{0.48\textwidth}
         \centering
         \includegraphics[width=\textwidth]{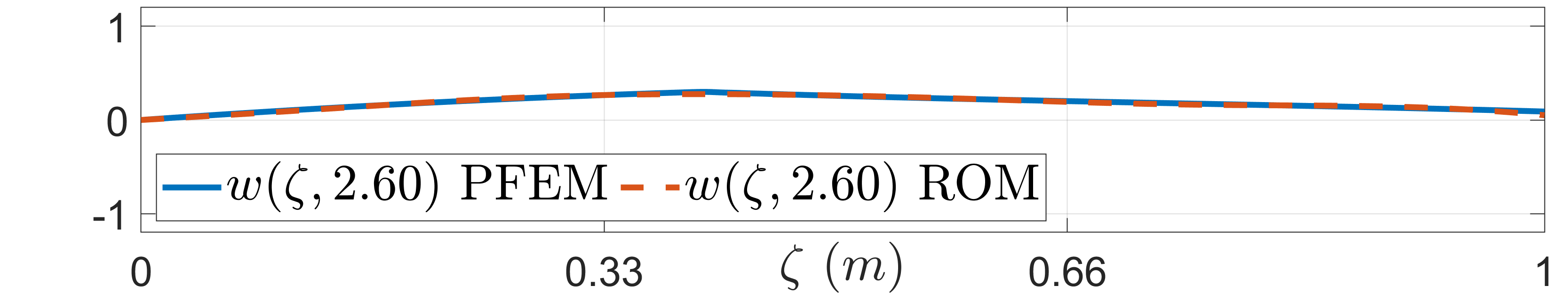}
%         \caption{}
%         \label{fig:SpectralZerosWave_b}
     \end{subfigure}
        \caption{Beam deformation of the discretized model and the ROM at different time steps.}
        \label{fig:SimulationBeam}
\end{figure}

\end{example}
\section{Conclusions and future work}\label{sec:Conclusion}

We have presented a systematic methodology for the discretization and reduction of one-dimensional BC-PHSs with first-order spatial derivatives. We cover a wide class of hyperbolic PDEs with a large class of boundary inputs and outputs. This is, for instance, the case of waves and beams with Neumann, Dirichlet, or mixed boundary conditions. We have proposed a structure preserving discretization scheme to derive a model that remains {\it passive} (or {\it impedance energy preserving}) if the initial system is. Finally, we have recalled the Loewner framework to find a ROM that interpolates the transfer function of the discretized model in a desired range of frequencies, obtaining similar input/output behaviour than the FOM at those frequencies. Moreover, to guarantee the passivity of the ROM, the interpolation algorithm can be repeated using the {\it spectral zeros} instead of pure frequency domain data. Additionally to this algorithm, we provide a constructive way to build a rectangular projection matrix that allows to approximate the states of the FOM from the states of the ROM, recovering the physical meaning of the state variables. 

Future works include the extension of the proposed approach to one-dimensional BC-PHSs with higher order spatial derivatives, including PDEs as the Euler-Bernoulli beam, which it has been already investigated for a particular case of boundary conditions in \citep{Cardoso2020JournalPartitioned}. Concerning the MOR of port-Hamiltonian systems, the case of Multiple-Input Multiple-Output (MIMO) can be investigated without using the tangential direction. Finally, the structure preserving MOR for parametric models is also part of the current work of the authors.



%%       \citet{<label>} ==> Jones et al. [21]
%%       \citep{<label>} ==> [21]
%%



\bibliographystyle{elsarticle-num-names} 
\bibliography{References}

\appendix

\small
\section{Weak Formulation}\label{Appendix_WF}
\subsection{Weak formulation}
Let denote $v_1(\zeta) \in H^1([a,b];\Real^{m})$ and $v_2(\zeta) \in H^1([a,b];\Real^{m})$ any arbitrary column vectors functions of sizes $m$. From \eqref{PDE_and_Constitutive}, we can obtain the following weak formulation multiplying by $v_1^\top$ from the right to the left the first PDE and by $v_2^\top$ the second one. Integrating both over the spatial domain we obtain:
\begin{equation}\label{WeakForm}
\begin{split}
\int_a^b v_1 ^\top \dot{x}_1 d\zeta &= \int_a^b v_1^\top \left[ {P}_{11} \partial_\zeta e_1 +{P}_{12} \partial_\zeta e_2
- {G}_{11}  e_1- {G}_{12} 
e_2 \right]
 d\zeta , \\
 \int_a^b v_2 ^\top \dot{x}_2 d\zeta &= \int_a^b v_2^\top \left[ {P}_{21} \partial_\zeta e_1 +{P}_{22} \partial_\zeta e_2
- {G}_{21}  e_1- {G}_{22} 
e_2 \right]
 d\zeta , \\
\int_a^b v_1^\top e_1 d\zeta &= \int_a^b v_1^\top \mathcal{H}_1 x_1 d\zeta, \\
\int_a^b v_2^\top e_2 d\zeta &= \int_a^b v_2^\top \mathcal{H}_2 x_2 d\zeta.
\end{split}
\end{equation}
Integrating by parts the first two terms of the right-hand side of the first two equations, we obtain:
\begin{equation}\label{WeakForm2}
\begin{split}
\int_a^b v_1 ^\top \dot{x}_1d\zeta &= \left[ v_1^\top {P}_{11}  e_1 + v_1^\top {P}_{12}  e_2 \right]_a^b 
- \int_a^b v_1^\top \left[ {G}_{11} e_1 + G_{12} e_2\right] d\zeta
\\& \quad - \int_a^b \left( \partial_\zeta v_1 \right) ^\top {P}_{11}  e_1 d\zeta - \int_a^b \left( \partial_\zeta v_1 \right) ^\top {P}_{12}  e_2 d\zeta, \\
\int_a^b v_2 ^\top \dot{x}_2 d\zeta &= \left[ v_2^\top {P}_{21}  e_1 + v_2^\top {P}_{22}  e_2 \right]_a^b 
- \int_a^b v_2^\top \left[ {G}_{21} e_1 + G_{22} e_2\right] d\zeta
\\& \quad - \int_a^b \left( \partial_\zeta v_2 \right) ^\top {P}_{21}  e_1 d\zeta - \int_a^b \left( \partial_\zeta v_2 \right) ^\top {P}_{22}  e_2 d\zeta. \\
\end{split}
\end{equation}
In particular, if $v_1 = e_1$ and $v_2 = e_2$, the previous equation with \eqref{InputOutput} satisfy \eqref{Balance}.
% become:
%\begin{equation}\label{Eq:WeakForm}
%\begin{split}
%\int_a^b e_1 ^\top \dot{x}_1d\zeta &= \left[ e_1^\top \overline{P}_1  e_2 \right]_a^b- \int_a^b \partial_\zeta e_1^\top \overline{P}_1  e_2 d\zeta\\& \quad + \int_a^b e_1^\top \overline{P}_0 e_2d\zeta, \\
%\int_a^b e_2^\top \dot{x}_2 d\zeta &= \left[ e_2^\top  \overline{P}_1  e_1 \right]_a^b - \int_a^b \partial_\zeta e_2^\top  \overline{P}_1  e_1d\zeta \\ & \quad+ \int_a^b e_2^\top  \overline{P}_0 e_1 d\zeta  - \int_a^b e_2^\top  \overline{G}_0 e_2d\zeta,
%\end{split}
%\end{equation}
%satisfying the same Hamiltonian balance \eqref{Balance} with inputs and outputs as in \eqref{InputOutput}. 
Indeed, the balance of the Hamiltonian \eqref{Eq:Hamiltonian} is obtained as follows:
\[
\begin{split}
\dot{H}(t) &= \int_a^b e ^\top \dot{x}d \zeta, \\
&= \int_a^b v_1 ^\top \dot{x}_1d \zeta + \int_a^b v_2 ^\top \dot{x}_2 d \zeta,\\
& = \left[ v^\top  {P}  e \right]_a^b - \int_a^b \left(\partial_\zeta v\right)^\top  {P} e d\zeta - \int_a^b v^\top  {G} e d\zeta, \\
& = \left[ e^\top  {P}  e \right]_a^b - \int_a^b \left(\partial_\zeta e\right)^\top  {P} e d\zeta - \int_a^b e^\top  {G} e d\zeta, \\
& =\dfrac{1}{2} \left[ e^\top  {P}  e \right]_a^b - \int_a^b e ^\top  {G} e d\zeta, \\
& =\dfrac{1}{2} \begin{bmatrix}
e(b)\\
e(a)
\end{bmatrix}^\top \begin{bmatrix}
P & 0 \\
0 & -P
\end{bmatrix}\begin{bmatrix}
e(b)\\
e(a)
\end{bmatrix} -\int_a^b e^\top  {G} e d\zeta, \\
& = \begin{bmatrix}
e(b)\\
e(a)
\end{bmatrix}^\top 
V_\mathcal{C}^\top V_\mathcal{B} \begin{bmatrix}
e(b)\\
e(a)
\end{bmatrix} -\int_a^b e^\top {G} e d\zeta, \\
& = y(t)^\top u(t) -\int_a^b e^\top  {G} e d\zeta \\
& \leq y(t)^\top u(t).
\end{split}
\]
where we have used \eqref{WeakForm2} with $v_1 = e_1$ and $v_2 = e_2$, then used %the properties $\overline{P}_1 = \overline{P}_1^\top$ and $\overline{P}_0 = -\overline{P}_0^\top$, then the fact that $\mathcal{R}^\top \Sigma \mathcal{R} = \left( \begin{smallmatrix} P_1 & 0_n \\ 0 _n & -P_1 \end{smallmatrix}\right)$, and finally 
Lemma~\ref{Lemma:InputOutput} with the fact that $\tilde{J}_\xi = -\tilde{J}_\xi^\top$. One can see that since $G+G^\top \leq 0$, the same balance equation \eqref{Balance} is obtained using the weak formulation.  

\end{document}